\input ./style/arxiv-general.cfg
\documentclass[aos,MSNbibl,seceqn,nameyear,dvips]{arximspdf}
\makeatletter
   \@ifpackageloaded{graphicx}{}{\usepackage{graphicx}}
\makeatother

\doi{10.1214/14-AOS1270}
\volume{43}
\issue{4}
\pubyear{2015}
\firstpage{1391}
\lastpage{1428}
\docsubty{FLA}
\relateddois{T11}{Discussed in
\relateddoi{d}{10.1214/14-AOS1270A},
\relateddoi{d}{10.1214/15-AOS1270B},
\relateddoi{d}{10.1214/15-AOS1270C},
\relateddoi{d}{10.1214/15-AOS1270D} and
\relateddoi{d}{10.1214/15-AOS1270E};
rejoinder at \relateddoi{r}{10.1214/15-AOS1270REJ}.}

\makeatletter

\def\vfrac#1#2{(#1)/#2}
\def\afrac#1#2{#1/(#2)}
\def\vafrac#1#2{(#1)/(#2)}

\newcommand{\rrvert}{\vert}
\newcommand{\rrVert}{\Vert}
\newcommand{\llvert}{\vert}
\newcommand{\llVert}{\Vert}
\def\diam{{\operatorname{diam}_n}}
\def\g{\gamma}
\def\Aupper{A}
\renewcommand{\b}{\beta}
\renewcommand{\k}{\kappa}
\newcommand{\EE}{{\mathbb{E}}}
\newcommand{\E}{\mathrm{E}}
\newcommand{\Var}{\operatorname{var}}
\newcommand{\sd}{\operatorname{sd}}
\newcommand{\eps}{\varepsilon}
\newcommand{\RR}{\mathbb{R}}
\newcommand{\NN}{\mathbb{N}}
\newcommand{\MM}{\mathbb{M}}
\renewcommand{\th}{\theta}
\newcommand{\s}{\sigma}
\renewcommand{\a}{\alpha}
\newcommand{\argmax}{\mathop{\operatorname{argmax}}}
\newtheorem{theorem}{Theorem}[section]
\newproclaim{definition}[theorem]{Definition}
\newtheorem{proposition}[theorem]{Proposition}
\newtheorem{lemma}[theorem]{Lemma}
\newproclaim{remark}[theorem]{Remark}
\makeatother

\begin{document}
\begin{frontmatter}

\title{Frequentist coverage of adaptive nonparametric Bayesian
credible sets\thanksref{T11,T2}}
\runtitle{Coverage of Bayesian credible sets}

\begin{aug}
\author[A]{\fnms{Botond}~\snm{Szab\'o}\ead[label=e1]{b.szabo@tue.nl}},
\author[B]{\fnms{A.~W.}~\snm{van der Vaart}\corref{}\ead[label=e2]{avdvaart@math.leidenuniv.nl}}
\and
\author[C]{\fnms{J. H.}~\snm{van Zanten}\ead[label=e3]{hvzanten@uva.nl}}
\runauthor{B. Szab\'o, A. W. van der Vaart and J. H. van Zanten}
\affiliation{TU Eindhoven, Leiden University and University of Amsterdam}
\address[A]{B. Szab\'o\\
Department of Mathematics\\
Eindhoven University of Technology\\
P.O. Box 513\\
5600 MB Eindhoven\\
The Netherlands\\
\printead{e1}}
\address[B]{A. van der Vaart\\
Mathematical Institute\\
Leiden University\\
P.O. Box 9512\\
2300 RA Leiden\\
The Netherlands\\
\printead{e2}}
\address[C]{H. van Zanten\\
Korteweg-de Vries Institute for Mathematics\\
University of Amsterdam\\
P.O. Box 94248\\
1090 GE Amsterdam\\
The Netherlands\\
\printead{e3}}
\end{aug}
\thankstext{T2}{Supported in part by the Netherlands Organization for
Scientific Research (NWO).
The research leading to these results has received funding from the
European Research Council
under ERC Grant Agreement 320637.}

%
\received{\smonth{11} \syear{2013}}
%
\revised{\smonth{4} \syear{2014}}

%
\begin{abstract}
We investigate the frequentist coverage of Bayesian credible sets
in a nonparametric setting. We consider a scale of priors of varying
regularity and choose the regularity by an empirical Bayes method.
Next we consider a central set of prescribed posterior probability
in the posterior distribution of the chosen regularity.
We show that such an adaptive Bayes credible set gives correct
uncertainty quantification of ``polished tail'' parameters,
in the sense of high probability of coverage of such parameters. On the negative
side, we show by theory and example that adaptation of the prior
necessarily leads to gross and haphazard uncertainty quantification for
some true parameters that are still within the hyperrectangle
regularity scale.
\end{abstract}

%
\end{frontmatter}

\section{Introduction}\label{sec1}
In Bayesian nonparametrics posterior distributions for functional
parameters are often
visualized by plotting a center of the posterior distribution, for
instance, the posterior mean
or mode, together with upper and lower bounds indicating a \emph
{credible set}, that is, a set that
contains a large fraction of the posterior mass (typically 95\%). The credible
bounds are intended to visualize the remaining uncertainty in the
estimate. In this paper
we study the validity of such bounds from a frequentist perspective in
the case of priors
that are made to adapt to unknown regularity.

It is well known that in infinite-dimensional models Bayesian credible
sets are not automatically
frequentist confidence sets, in the sense that under the assumption
that the data are in actual fact
generated by a ``true parameter,'' it is not automatically true that
they contain that truth with
probability at least the credible level.
The earliest literature focused on negative examples, showing that
for many combinations of truths and priors, Bayesian credible sets can
have very bad or at least
misleading frequentist behavior; see, for instance, \citet{Cox},
\citet
{Freedman},
\citet{Johnstone}. [An exception is \citet{Wahba}, who showed
encouraging simulation
results and gives heuristic arguments for good performance.]
However, credible sets do not always have bad frequentist coverage. In
the papers
\citet{Leahu}, \citeauthor{Bartek} (\citeyear{Bartek,Bartek2}) this
matter was
investigated in the setting of
the canonical (inverse) signal-in-white-noise model, where,
essentially, the unknown parameter was a
function with a fixed regularity and the (Gaussian) prior had a fixed
regularity as
well. The main message in these papers is that Bayesian credible sets
typically have {\sl good}
frequentist coverage in case of undersmoothing (using a prior, i.e.,
less regular
than the truth), but can have coverage zero and be far too small in the
other case. Simulation studies
corroborate these theoretical findings and show that the problem of misleading
uncertainty quantification is a very practical one.

The solution to undersmooth the truth, which gives good uncertainty
quantification, is unattractive for
two reasons. First, it leads to a loss in the quality of the
reconstruction, for example, by the posterior
mode or mean. Second, the true regularity of the functional parameter
is never known and hence
cannot be used to select a prior that undersmoothes the right
regularity. Therefore, in practice, it is common to try
and ``estimate'' the regularity from the data, and thus to {\sl adapt}
the method
to the unknown regularity. Bayesian versions of this approach
can be implemented using empirical or hierarchical Bayes methods.
Empirical Bayes methods
estimate the unknown regularity using the marginal likelihood for the
data in the Bayesian
setup; see Section~\ref{SectionEB} for a precise description.
Hierarchical Bayes methods
equip the regularity parameter with a prior and follow a full Bayesian approach.

In the present paper we concentrate on the empirical Bayes approach.
In the context of the inverse signal-in-white-noise model, this method has
been shown to be rate-adaptive, in the sense that the posterior contracts
at a (near) optimal rate around the truth
for a range of true regularities, without using information about this
regularity [see \citet{KSzVZ}
for an analysis of the method in the present paper
or \citet{Ray} for similar work]. However, these papers only address
contraction of the posterior and do not
investigate frequentist coverage of credible sets, which is perhaps
more important
than rate-adaptiveness to validate the use of these methods. In the present
paper we study whether the
empirical Bayes method, which is optimal from the point of view of contraction
rates, also performs well from the perspective of coverage. In
particular, we investigate to which
extent the method yields adaptive confidence sets.

Bayesian credible sets can of course not beat the general fundamental
limitations of adaptive
confidence sets. As pointed out by \citet{Low}, it is in general not
possible to construct confidence
sets that achieve good coverage across a range of nested models with
varying regularities and at the
same time possess a size of optimal order when the truth is assumed to
be in one of the particular
sub-models. Similar statements, in various contexts, can be found in
\citet{JudLam},
\citeauthor{CaiLow} (\citeyear{CaiLow,CaiLow2006}), \citet{RobVaart}, \citet{GenWas}, and \citet{Hoff}.

We show in this paper that for the standard empirical Bayes procedure
(which is rate-adaptive), there
always exist truths that are not covered asymptotically by its credible
sets. This bad news
is alleviated by the fact that there are only a few of these
``inconvenient truths'' in some sense. For
instance, the minimax rate of estimation does not improve after
removing them from the model;
they form a small set in an appropriate topological sense; and they
are unlikely under any of the priors. The good news is that after removing
these bad truths, the empirical Bayes credible sets become adaptive
confidence sets with good
coverage.

Our results are inspired by recent (non-Bayesian) results of
\citet{GineNickl} and \citet{Bull}. These authors also
remove a ``small'' set of undesirable truths from the model and focus
on so-called self-similar
truths. Whereas these papers use theoretical frequentist methods of adaptation,
in the present paper our starting point is the Bayesian (rate-adaptive)
procedure. This generates candidate confidence sets
for the true parameter (the credible sets), that are routinely used in practice.
We next ask for which truths this practice can be justified and for
which not.
Self-similar truths, defined appropriately in our setup, are covered,
but also a more general class of parameters, which we call \emph
{polished tail sequences}.

The paper is structured as follows. In the next section we describe the setting:
the inverse signal-in-white-noise model and the adaptive empirical
Bayes procedure. In Section~\ref{secmain} the associated credible sets
are constructed and analyzed. A
first theorem exhibits truths
that are not covered asymptotically by these sets.
The second theorem shows that when these ``inconvenient
truths'' are removed, the credible sets yield adaptive, honest
confidence sets. The theoretical
results are illustrated in Section~\ref{secsim} by a simulation
study. Proofs are given in
Sections~\ref{secproofcov}--\ref{SectionTechnicalLemmas}.
In Section~\ref{seccon} we conclude with some
remarks about possible extensions and generalizations.
Finally, the \hyperref[secproofMLE2]{Appendix} is a self-contained
proof of a version of an important auxiliary result first proved in
\citet{KSzVZ}.

We conclude the \hyperref[sec1]{Introduction} by further discussions of adaptive nonparametric
confidence sets and the coverage of credible sets.

The credible sets we consider in this paper are $\ell_2$-balls, even
though we believe
that similar conclusions will be true for sets of different shapes. It
is known that
$\ell_2$-confidence balls can be honest over a model of regularity $\a$
and possess a radius that adapts to the minimax rate whenever the true
parameter is
of smoothness contained in the interval $[\a, 2\a]$. Thus, these
balls can adapt
to double a coarsest smoothness level [\citet{JudLam}, \citet
{CaiLow2006}, \citet{RobVaart}, \citet{BullNickl}].
The fact that the coarsest level $\a$ and the radius of the ball must
be known [as shown in
\citet{BullNickl}] makes this type of adaptation somewhat theoretical.
This type of adaptation is not considered in the present paper (in
fact, we do not have
a coarsest regularity level $\a$).
Work we carried out subsequent to the present paper [\citet{SzaboJSPI}]
indicates that this type of adaptation
can be incorporated in the Bayesian framework, but requires a different
empirical Bayes procedure
as the one in the present paper [based on the likelihood (\ref
{EqDefinitionHata})].
Interestingly, with the latter method
and for known $\a$, adaptation occurs for all true parameters, also
the inconvenient ones.

The credible sets considered in the present paper result from posterior
distributions
for infinite-dimensional parameters, or functions, which implicitly
make the bias--variance
trade-off that is characteristic of nonparametric estimation. These
posterior distributions
induce marginal posterior distributions for real-valued functionals of
the parameter. If such a
functional is sufficiently smooth, then the corresponding marginal
posterior distribution may satisfy
a Bernstein--von Mises theorem, which typically entails that the bias is
negligible relative
to the variance of estimation. Just as in the case of
finite-dimensional models, such
an approximation implies that the credible sets for the functional are
asymptotically equivalent to frequentist confidence sets. In this sense
nonparametric priors and posteriors may yield exact, valid credible sets.
By extending this principle to an (infinite) collection of smooth
functionals that identifies
the parameter, \citet{CastilloNickl} even obtain an exact credible set
for the full parameter.
However, elegant as their construction may be,
it seems that no method that avoids dealing with the bias--variance trade-off
will properly quantify the uncertainty of nonparametric Bayesian
inference as it is
applied in current practice.

\subsection{Notation}\label{SectionNotation}
The $\ell^2$-norm of an element $\th\in\ell^2$ is denoted by $\llVert
\th
\rrVert $, that is, $\llVert \th\rrVert ^2=\sum_{i=1}^{\infty}\th
_i^2$. The \emph{hyperrectangle} and
\emph{Sobolev space} of order $\b>0$ and (square) radius $M > 0$ are
the sets
%
%
\begin{eqnarray}
\Theta^{\b}(M)&=&\Bigl\{\th\in\ell^2\dvtx  \sup
_{i\ge1}i^{1+2\b}\th_i^2 \le M
\Bigr\},\label{defhyperrectangle}
\\
S^{\b}(M)&=&\Biggl\{\th\in\ell^2\dvtx  \sum
_{i=1}^\infty i^{2\b}\th_i^2
\le M\Biggr\}.\label{defsobolevspace}
\end{eqnarray}
For two sequences $(a_n)$ and $(b_n)$ of numbers, $a_n
\asymp b_n$ means that $\llvert a_n/b_n\rrvert $ is bounded away from
zero and
infinity, $a_n \lesssim b_n$
that $a_n/b_n$ is bounded, $a_n \sim b_n$ that $a_n/b_n \rightarrow1$, and
$a_n \ll b_n$ that $a_n / b_n \rightarrow0$, all as $n$ tends to
infinity. The maximum and minimum
of two real numbers $a$ and $b$ are denoted by
$a \vee b$ and $a \wedge b$.

\section{Statistical model and adaptive empirical Bayes procedure}\label{SectionEB}
We formulate and prove our results in the canonical setting of the
inverse signal-in-white-noise
model. As usual, we reduce it to the sequence formulation. See, for
instance, \citet{Cavalier} and the
references therein for more background and many examples fitting this framework.

The observation is a sequence $X=(X_1,X_2,\ldots)$ satisfying
%
%
\begin{equation}
\label{eqdata} X_i=\k_i\th_{0,i}+
\frac{1}{\sqrt{n}}Z_i, \qquad i = 1, 2, \ldots,
\end{equation}
where $\th_{0}=(\th_{0,1},\th_{0,2},\ldots)\in\ell^2$ is the unknown
parameter of interest,
the $\k_i$'s are known constants (transforming the truth) and the
$Z_i$ are independent, standard normally distributed random variables.
The rate of decay of the $\k_i$'s determines the difficulty of the
statistical problem of recovering $\th_0$. We
consider the so-called mildly ill-posed case where
%
%
\begin{equation}
C^{-2}i^{-2p}\leq\k_i^2\leq
C^2 i^{-2p}\label{assumpinverse},
\end{equation}
for some fixed $p \ge0$ and $C > 0$. In particular, the choice $p=0$
corresponds to the
ordinary signal-in-white-noise model, whereas $p>0$ gives a true \emph
{inverse problem}.

For $\a> 0$ we define a prior measure $\Pi_\a$ for the parameter
$\th_0$ in (\ref{eqdata}) by
%
%
\begin{equation}
\label{prior} \Pi_{\a}=\bigotimes_{i=1}^{\infty}N
\bigl(0,i^{-1-2\a}\bigr).
\end{equation}
The coordinates $\th_i$ are independent under this prior. Since the
corresponding
coordinates of the data are also independent, the independence is retained
in the posterior distribution, which by univariate conjugate Gaussian
calculation
can be seen to be
%
%
\begin{equation}
\label{PostDist} \Pi_\a(\cdot\mid X) = \bigotimes
_{i=1}^{\infty}{N} \biggl(\frac
{n\k_i^{-1}}{i^{1+2\a}\k_i^{-2} + n}X_i,
\frac{\k_i^{-2}}{i^{1+2\a}\k_i^{-2}+n} \biggr).
\end{equation}
The prior (\ref{prior}) puts mass 1 on Sobolev spaces and
hyperrectangles of every order strictly
smaller than $\a$ (see Section~\ref{SectionNotation} for
definitions), and hence expresses a prior
belief that the parameter is regular of order (approximately) $\a$.
Indeed, it is shown in
\citet{Bartek} that if the true parameter $\th_0$ in (\ref{eqdata})
belongs to a Sobolev space of
order $\a$, then the posterior distribution contracts to the true
parameter at the minimax rate
$n^{-2\a/(1+2\a+2p)}$ for this Sobolev space. A similar result can be
obtained for
hyperrectangles. On the other hand, if the regularity of the true
parameter is different from $\a$,
then the contraction can be much slower than the minimax rate.

The suboptimality in the case the true regularity is unknown
can be overcome by a data-driven choice of $\a$. The
\emph{empirical Bayes} procedure consists in replacing the fixed
regularity $\a$ in
(\ref{PostDist}) by (for given $\Aupper$, possibly dependent on $n$)
%
%
\begin{equation}
\label{EqDefinitionHata} \hat\a_n=\argmax_{\a\in[0,\Aupper]}
\ell_n(\a),
\end{equation}
where $\ell_n$ is the marginal log-likelihood for $\a$ in the
Bayesian setting:
$\th\mid\a\sim\Pi_\a$ and $X\mid(\th, \a)\sim\bigotimes_{i}{N}(\k_i\th_i,1/n)$.
This is given by
%
%
\begin{equation}
\ell_n(\a)=-\frac{1}{2}\sum_{i=1}^{\infty}
\biggl( \log\biggl(1+\frac{n}{i^{1+2\a}\k_i^{-2}} \biggr)-\frac{n^2}
{i^{1+2\a}\k
_i^{-2}+n}X_i^2
\biggr).\label{fnalpha}
\end{equation}
If there exist multiple maxima, any one of them can be chosen.

The \emph{empirical Bayes posterior} is defined as the random measure
$\Pi_{\hat\a_n}(\cdot\mid X)$
obtained by substituting $\hat\a_n$ for $\a$ in the posterior
distribution (\ref{PostDist}), that is,
%
%
\begin{equation}
\label{EqEBPosteriorDistribution} \Pi_{\hat\a_n}(\cdot\mid X) = \Pi_{\a
}(\cdot\mid
X) |_{\a= \hat\a_n}.
\end{equation}
In \citet{KSzVZ} this distribution is shown
to contract to the true parameter at the (near) minimax rate within the
setting of Sobolev
balls and also at the near optimal rate in situations of supersmooth parameters.
Extension of their results shows that the posterior distribution
also performs well for many other models, including hyperrectangles.
Thus, the empirical Bayes posterior distribution manages to
recover the true parameter by \emph{adapting} to unknown models.

We now turn to the main question of the paper: can the spread
of the empirical Bayes posterior distribution be
used as a measure of the remaining uncertainty in this recovery?

\subsection{Notational assumption}
We shall from now on assume that the first coordinate $\th_{0,1}$ of the
parameter $\th_0$ is zero. Because the prior (\ref{prior}) induces a
$N(0,1)$ prior on $\th_{0,1}$,
which is independent of $\a$, the marginal likelihood function (\ref
{fnalpha}) depends on $X_1$
only through a vertical shift, independent of $\a$. Consequently, the
estimator $\hat\a_n$ does
not take the value of $\th_{0,1}$ into account. While this did not
cause problems for the minimax
adaptivity mentioned previously, this does hamper the performance of
credible sets obtained from the
empirical Bayes posterior distribution, regarding uniformity in the
parameter. [In fact, the posterior distribution for $\th_1$ has mean
equal to
$n\kappa_1^{-1}/(\kappa_1^{-2}+n)X_1$, independent of $\a$. The bias
of this estimator
of $\th_{0,1}$ for fixed $n$ would lead to arbitrarily small coverage
for values
$\th_{0,1}\rightarrow\pm\infty$, invalidating the main result of
the paper, Theorem~\ref{thmCoverage} below.]
One solution would be to use the variances $(i+1)^{-1-2\a}$
in (\ref{prior}). For notational simplicity we shall instead assume
that $\th_{0,1}=0$ throughout the remainder of the paper.

\section{Main results: Asymptotic behavior of credible sets}\label{secmain}
For fixed $\a> 0$, let $\hat\th_{n, \a}$ be the posterior mean
corresponding to the prior $\Pi_\a$
[see (\ref{PostDist})]. The centered posterior is a Gaussian measure
that does not depend on the
data and, hence, for $\gamma\in(0,1)$ there exists a deterministic
radius $r_{n,\gamma}(\a)$ such
that the ball around the posterior mean with this radius receives a
fraction $1-\gamma$ of the
posterior mass, that is, for $\a>0$,
%
%
\begin{equation}
\label{eqradius} \Pi_\a\bigl(\th\dvtx  \llVert\th- \hat
\th_{n, \a}\rrVert\le r_{n, \gamma}(\a) \mid X \bigr) = 1-\gamma.
\end{equation}
In the exceptional case that $\a=0$, we define the radius to be infinite.
The empirical Bayes credible sets that we consider in this paper are
the sets obtained by replacing the fixed
regularity $\a$ by the data-driven choice $\hat\a_n$. Here we
introduce some more flexibility
by allowing the possibility of blowing up the balls by a factor $L$.
For $L > 0$ we define
%
%
\begin{equation}
\label{EqCredibleSet} \hat C_n(L) = \bigl\{\th\in\ell^2\dvtx
\llVert\th- \hat\th_{n, \hat\a
_n}\rrVert\le Lr_{n, \gamma}(\hat
\a_n) \bigr\}.
\end{equation}
By construction, $\Pi_{\hat\a_n}(\hat C_n(L) \mid X) \ge1-\gamma
$ iff $L \ge1$.

We are interested in the performance of the random sets $\hat C_n(L)$
as frequentist confidence
sets. Ideally, we would like them to be \emph{honest} in the sense that
\[
\inf_{\th_0\in\Theta_0}P_{\th_0} \bigl(\th_0 \in\hat
C_n(L) \bigr) \ge1-\gamma,
\]
for a model $\Theta_0$ that contains all parameters deemed possible.
In particular, this model should contain parameters of all regularity levels.
At the same time we would like the sets to be \emph{adaptive},
in the sense that the radius of $\hat C_n(L)$ is
(nearly) bounded by the optimal rate for a model of a given regularity level,
whenever $\th_0$ belongs to this model.
As pointed out in the \hyperref[sec1]{Introduction}, this is too much to ask, as
confidence sets with
this property, Bayesian or non-Bayesian, do not exist. For the present
procedure we can explicitly
exhibit examples of ``inconvenient truths'' that are not covered at all.

%
\begin{theorem}\label{LemmaCounterExample}
For given $\beta, M>0$ and $1\le\rho_j\uparrow\infty$, and
positive integers $n_j$
with $n_{j+1}\geq(2\rho_{j+1}^2)^{1+2\beta+2p}n_j$,
define $\th_0=(\th_{0,1},\th_{0,2},\ldots)$ by
\begin{eqnarray*}
\th_{0,i}^2=\cases{ 0, &\quad if $\rho_j^{-1}n_j^{\afrac{1}{1+2\b+2p}}
\leq i<n_j^{\afrac
{1}{1+2\b+2p}}$, $j=1,2,\ldots,$
\vspace*{3pt}\cr
0, &\quad if
$2n_j^{\afrac{1}{1+2\b+2p}}\leq i\leq\rho_jn_j^{\afrac
{1}{1+2\b+2p}}$,
$j=1,2,\ldots,$
\vspace*{3pt}\cr
Mi^{-1-2\beta}, &\quad otherwise.}
\end{eqnarray*}
Then the constant $M$ can be chosen such that $P_{\th_0}(\th_0 \in
\hat C_{n_j}(L_{n_j})) \rightarrow0$ as
$j \rightarrow\infty$ for every $L_{n_j}\lesssim\sqrt{M} \rho
_j^{(1+2p)/(8+8\beta+8p)}$.
\end{theorem}

For the proof see Section~\ref{secCounterExample}.

By construction, the (fixed) parameter $\th_0$ defined in Theorem~\ref{LemmaCounterExample}
belongs to the hyperrectangle $\Theta^\b(M)$,
and in this sense is ``good,'' because of the ``smooth'' truth.
However, it is an
inconvenient truth, as it tricks the empirical Bayes procedure, making
this choose the ``wrong''
regularity $\a$, for which the corresponding credible set does not
cover $\th_0$. The intuition
behind this counterexample is that for a given sample size or noise
level $n$ the empirical Bayes
procedure, and any other statistical method, is able to judge the coordinates
$\th_{0,1},\th_{0,2},\ldots$ only up to a certain effective
dimension $N_n$, fluctuations in the
higher coordinates being equally likely due to noise as to a nonzero
signal. Now if
$(\th_{0,1},\ldots, \th_{0,N_n})$ does not resemble the infinite sequence
$(\th_{0,1},\th_{0,2},\ldots)$, then the empirical Bayes procedure
will be tricked into choosing a
smoothness $\hat\a_n$ that does not reflect the smoothness of the
full sequence, and failure of
coverage results. The particular example $\th_0$ in
Theorem~\ref{LemmaCounterExample} creates this situation by including
``gaps'' of 0-coordinates. If the effective dimension is at the end of
a gap of such 0-coordinates,
then the empirical Bayes procedure will conclude that $\th_0$ is
smoother than it really is, and
make the credible set too narrow.

A more technical explanation can be given in terms of a bias--variance
trade-off. For the given truth in Theorem~\ref{LemmaCounterExample}
the bias of the posterior mean
is of ``correct order'' $n^{-\b/(1+2\b+ 2p)}$ corresponding to
$\Theta^\b(M)$,
but (along the subsequence $n_j$) the spread of the posterior is of
strictly smaller
order. Thus, the posterior distribution is overconfident about its performance.
See Section~\ref{secCounterExample} for details.

Intuitively, it is not surprising that such bad behavior occurs, as
nonparametric
credible or confidence sets always necessarily extrapolate into aspects
of the truth
that are not visible in the data. Honest uncertainty quantification is
only possible by
a priori assumptions on those latter aspects. In the context of
regularity, this
may be achieved by ``undersmoothing,'' for instance, by using a prior
of fixed regularity smaller than the true regularity. Alternatively, we
may change the notion
of regularity and strive for honesty over different models. In the latter
spirit we shall show that the empirical Bayes credible sets $\hat C_n(L)$
are honest over classes of ``polished'' truths.

%
\begin{definition}
A parameter $\th\in\ell^2$ satisfies the \emph{polished tail
condition} if,
for fixed positive constants $L_0$, $N_0$ and $\rho\ge2$,
%
%
\begin{equation}
\label{EqPolishedTail} \sum_{i=N}^{\infty}
\th_i^2\le L_0 \sum
_{i=N}^{\rho N}\th_i^2\qquad \forall N\ge N_0.
\end{equation}
\end{definition}

We denote by $\Theta_{pt}(L_0,N_0,\rho)$ the set of all
polished tail sequences $\th\in\ell^2$
for the given constants $L_0$, $N_0$ and $\rho$. As the constants
$N_0$ and $\rho$ are fixed in most of the following (e.g., at $N_0=2$
and $\rho=2$),
we also use the shorter $\Theta_{pt}(L_0)$. [It would be possible to
make a refined
study with $N_0$ and $L_0$ tending to infinity, e.g., at
logarithmic rates, in order
to cover a bigger set of parameters, eventually. The result below would
then go through
provided the constant $L$ in the credible sets $C_n(L)$ would also tend
to infinity at a related rate.]

The condition requires that the contributions of the blocks $(\th
_N,\ldots,\th_{N\rho})$
of coordinates to the $\ell^2$-norm of $\th$ cannot surge over the
contributions
of earlier blocks as $N\rightarrow\infty$.
Sequences $\th$ of exact polynomial order $\th_{i}^2\asymp i^{-1-2\b
}$ on the
``boundary'' of hyperrectangles are obvious examples,
but so are the sequences $\th_i\asymp i^qe^{-\zeta i^c}$ and the sequences
$\th_{i}^2\asymp(\log i)^q i^{-1-2\b}$ (for $q\in\RR$, $\zeta,c>0$).
Furthermore, because the condition is not on the
individual coordinates $\th_i$, but on blocks of coordinates of
increasing size,
the set of polished tail sequences is in fact much larger than these
coordinatewise regular examples suggest.

In particular, the set includes the ``self-similar'' sequences. These
were defined
by \citet{PicTri} and employed by \citet{GineNickl} and
\citet{Bull},
in the context of
wavelet bases and uniform norms. An $\ell^2$ definition in
the same spirit with reference to hyperrectangles is as follows.

%
\begin{definition}
A parameter $\th\in\Theta^{\b}(M)$ is \emph{self-similar} if,
for some fixed positive constants $\eps$, $N_0$ and $\rho\ge2$,
%
%
\begin{equation}
\sum_{i=N}^{\rho N}\th_{i}^2
\geq\eps M N^{-2\b}\qquad\forall N\ge N_0. \label{propSelfSim}
\end{equation}
\end{definition}

We denote the class of self-similar elements of $\Theta^\b(M)$ by
$\Theta_{ss}^\b(M,\eps)$. The
parameters $N_0$ and $\rho$ are fixed and omitted from the notation.

If we think of $\th$ as a sequence of Fourier coefficients, then the
right-hand side of (\ref{propSelfSim}) without $\eps$ is the maximal
energy at frequency $N$ of a sequence in a hyperrectangle
of radius $\sqrt M$. Thus, (\ref{propSelfSim}) requires that the
total energy in
every \emph{block} of consecutive frequency components is a fraction
of the energy of a typical
signal: the signal looks similar at all frequency levels. Here the
blocks increase
with frequency (with lengths proportional to frequency), whence the
required similarity is only
on average over large blocks.

Self-similar sequences are clearly polished tail sequences, with
$L_0=\eps^{-1}$.
Whereas the first refer to a particular regularity class, the latter do
not. As polished tail sequences
are defined by self-referencing, they might perhaps be considered
``self-similar'' in
a generalized sense. We show in Theorem~\ref{thmCoverage}
that the polished tail condition is sufficient
for coverage by the credible sets (\ref{EqCredibleSet}).
Self-similarity is restrictive.
For instance, the polished tail sequence\vspace*{1pt} $\th_i=i^{-\b-1/2}(\log
i)^{-q/2}$ is contained in the hyperrectangle
$\Theta^\b(1)$ for every $q\ge0$, and also in $S^\b(M)$ for some
$M$ if $q>1$, but
it is not self-similar for any $q>0$. This could be remedied by
introducing (many) different
types of self-similar sequences, but the self-referencing of polished
tail sequences seems
much more elegant.

%
\begin{remark}
An alternative to condition (\ref{propSelfSim}) would be the
slightly weaker
\[
\sum_{i=N}^{\infty}\theta_i^2
\geq\eps M N^{-2\beta}\qquad\forall N\ge N_0.
\]
This removes the parameter $\rho$, but, as $\theta$ is assumed to be contained
in the hyperrectangle $\Theta^\b(M)$, it can be seen that this
seemingly relaxed condition
implies (\ref{propSelfSim}) with $\eps$ replaced by $\eps/2$ and
$\rho$ sufficiently large that
$\sum_{i>\rho N} i^{-1-2\b}< \eps/2 N^{-2\b}$.
\end{remark}

One should ask how many parameters are not polished tail or self-similar.
We give three arguments that there are only few: topological,
minimax, and Bayesian.

A topological comparison of the classes of self-similar and
non self-similar functions
obviously depends on the chosen topology. From the proof of
Theorem~\ref{LemmaCounterExample}
it is clear that the lack of coverage is due to the tail behavior of
the non self-similar
(or nonpolished-tail) truth in the statement of the theorem. Hence, by
modifying the tail
behavior of an arbitrary sequence $\th$ we can get a non self-similar
sequence with asymptotic
coverage 0. Similarly, every element of $\Theta^\b(M)$ can be made
self-similar by modifying
its tail. So we see that in the $\ell^2$-norm topology the difference
in size between
the two classes does not become apparent. Both the self-similar and the
bad, non self-similar
truths are dense in $\Theta^\b(M)$. Following \citet{GineNickl}, one
can also consider matters
relative to the finer smoothness topology on $\Theta^\b(M)$. Similar
to their Proposition~4, it
can then be shown that the set of non self-similar functions is nowhere
dense, while the set of
self-similar functions is open and dense. This suggests that
self-similarity is the norm rather than
the exception, as is also expressed by the term \emph{generic} in the
topological sense.

The minimax argument for the neglibility of nonpolished-tail sequences
is that restriction
to polished tail (or self-similar) truths does not reduce the
statistical difficulty of the problem.
We show below (see Proposition~\ref{LemmaMinimaxRate}) that
restriction to self-similarity changes
only the constant in the minimax risk for hyperrectangles and not the
order of magnitude
or the dependence on the radius of the rectangle.
Similarly, the minimax risk over Sobolev balls is reduced by at most
a logarithmic factor by a restriction to polished tail sequences
(see Proposition~\ref{PropositionSobolevSize}).

A third type of reasoning is that polished tail sequences are \emph
{natural} once
one has adapted the Bayesian setup with priors of the form (\ref
{prior}). The following
proposition shows that almost every realization from such a prior is a
polished tail
sequence for some $N_0$. By making $N_0$ large enough
we can make the set of polished tail sequences have arbitrarily large
prior probability.
This is true for any of the priors $\Pi_\a$ under consideration.
Thus, if one believes
one of these priors, then one accepts the polished tail condition.
The result may be compared to Proposition~4 of \citet{Hoff} and
Proposition~2.3 of
\citet{Bull}.

Recall that $\Theta_{pt}(L_0,N_0,\rho)$ is the set of $\th\in\ell
_2$ that satisfy~(\ref{EqPolishedTail}).

%
\begin{proposition}
For every $\a>0$ the prior $\Pi_\a$ in (\ref{prior}) satisfies\break
$\Pi_\a(\bigcup_{N_0}\Theta_{pt}(2/\a+1,N_0,2) )=1$.
\end{proposition}

\begin{pf}
Let $\th_1,\th_2,\ldots$ be independent random variables with $\th
_i\sim N(0,\break  i^{-1-2\a})$,
and let\vspace*{1pt} $\Omega_N$ be the event $\{\sum_{i\ge N}\th_i^2>(2/\a
+1)\sum_{i=N}^{2N}\th_i^2\}$.
By the Borel--Cantelli lemma it suffices to show that $\sum_{N\in\NN}
\Pi_\a(\Omega_N)<\infty$.
We have that
\begin{eqnarray*}
\E\Biggl(\frac{2+\a}{\a}\sum_{i=N}^{2N}
\th_i^2-\sum_{i \ge N}\th
_i^2 \Biggr) &=&\frac{2+\a}{\a}\sum
_{i=N}^{2N} \frac{1}{i^{1+2\a}}-\sum
_{i\ge
N}\frac{1}{i^{1+2\a}}
\\
&\ge&\frac{2}{\a}\sum_{i=N}^{2N}
\frac{1}{i^{1+2\a}}-\int_{2N}^{\infty}x^{-1-2\a} \,dx
\\
&\ge&\bigl(2^{-2\a-1}/\a\bigr) N^{-2\a}.
\end{eqnarray*}
Therefore, by Markov's inequality, followed by the
Marcinkiewitz--Zygmund and H\"older
inequalities, for $q\ge2$ and $r>1$,
\begin{eqnarray*}
\Pi_\a(\Omega_N) &\lesssim& N^{2\a q}\E\Biggl
\llvert \sum_{i=N}^{2N}\bigl(
\th_i^2-\E\th_i^2\bigr) -\sum
_{i>2N}\bigl(\th_i^2-\E
\th_i^2\bigr)\Biggr\rrvert^q
\\
&\lesssim& N^{2\a q} \E\biggl(\sum_{i\ge N}
\bigl(\th_i^2-\E\th_i^2
\bigr)^2 \biggr)^{q/2}
\\
&\lesssim& N^{2\a q} \sum_{i\ge N} \E\bigl(
\th_i^2-\E\th_i^2
\bigr)^qi^{r(q/2-1)} \biggl(\sum_{i\ge N}i^{-r}
\biggr)^{q/2-1}.
\end{eqnarray*}
Since $\E(\th_i^2-\E\th_i^2)^q\asymp i^{-2q\a-q}$, for $-2q\a
-q+r(q/2-1)<-1$
(e.g., $q>2$ and $r$~close to 1), the\vspace*{1pt} right-hand side is of the
order $N^{2\a q} N^{-2q\a-q+r(q/2-1)+1} \times\break  N^{-(r-1)(q/2-1)}=N^{-q/2}$.
This is summable for $q>2$.
\end{pf}

The next theorem is the main result of the paper. It states that when
the parameter is restricted to
polished tail sequences, the empirical Bayes credible ball $\hat
C_n(L)$ is an honest,
frequentist confidence set, if $L$ is not too small. In
Sections~\ref{SectionHyperrectangle}--\ref{SectionSupersmooth}
this theorem will be complemented
by additional results to show that $\hat C_n(L)$
has radius $r_{n,\gamma}(\hat{\a}_n)$ of minimax order over a range
of regularity classes.

Recall that $\Theta_{pt}(L_0)$ is the set of all polished tail
sequences $\th\in\ell^2$
for the given constant $L_0$, and $\Aupper$ is the constant in
(\ref{EqDefinitionHata}).

%
\begin{theorem}\label{thmCoverage}
For any $\Aupper, L_0, N_0$ there exists a constant $L$ such that
%
%
\begin{equation}
\inf_{\th_0\in\Theta_{pt}(L_0)} P_{\th_0} \bigl(\th_0\in\hat
C_n(L) \bigr)\rightarrow1.\label{eqPostCov}
\end{equation}
Furthermore,\vspace*{1pt} for $\Aupper=\Aupper_n\leq\sqrt{\log n}/(4\sqrt{\log
\rho\vee e})$ this is true with
a slowly varying sequence [$L:=L_n\lesssim(3\rho^{3(1+2p)})^{\Aupper
_n}$ works].
\end{theorem}

\begin{pf}
See Section~\ref{secMain}.
\end{pf}

The theorem shows that the sets $\hat C_n(L)$ are large enough to
catch any truth that satisfies the polished tail condition, in the sense
of \emph{honest} confidence sets. With the choice $L$ as
in the theorem, their coverage in fact tends to 1, so that they are
\emph{conservative}
confidence sets. In \citet{Bartek} it was seen that for
deterministic choices
of $\a$, the constant $L$ cannot be adjusted so that exact coverage
$1-\g$ results;\vspace*{1pt}
and $L$ in (\ref{eqPostCov}) may have to be larger than 1.
Thus, the credible sets $\hat C_n(L)$
are not exact confidence sets, but in combination with the
results of Sections~\ref{SectionHyperrectangle}--\ref{SectionSupersmooth}
the theorem does indicate that their
order of magnitude is correct in terms of frequentist confidence statements.

As the credible sets result from a natural Bayesian procedure, this message
is of interest by itself. In the next subsections we complement this
by showing that the good coverage is not obtained
by making the sets $\hat C_n(L)$ unduly large. On the contrary, their radius
$Lr_{n,\g}(\hat\a_n)$ is of the minimax estimation rate for various
types of models.
In fact, it is immediate from the definition of the radius in (\ref
{eqradius}) that
$r_{n,\g}(\hat\a_n)=O_P(\eps_n)$, whenever the posterior
distribution contracts
to the true parameter at the rate $\eps_n$, in the sense that for
every $M_n\rightarrow\infty$,
\[
E_{\th_0} \Pi_{\hat\a_n} \bigl(\th\dvtx  \llVert\th-\th_0
\rrVert> M_n\eps_n\mid X \bigr) \rightarrow0.
\]
Such contraction was shown in \citet{KSzVZ} to take place
at the (near) minimax rate $\eps_n=n^{-\b/(2\b+2p+1)}$ uniformly in
parameters ranging
over Sobolev balls $S^\b(M)$,
adaptively in the regularity level $\b$. In the next subsections we
refine this in various ways:
we also consider other models, and give refined and oracle statements for
the behavior of the radius under polished tail or self-similar sequences.

Generally speaking, the size of the credible sets $\hat C_n(L)$ are
of (near) optimal size whenever the empirical
Bayes posterior distribution (\ref{EqEBPosteriorDistribution})
contracts at the (near) optimal rate.
This is true for many but not all possible models for two reasons.
On the one hand, the choice of priors with variances $i^{-1-2\a}$, for
some $\a$, is
linked to a certain type of regularity in the parameter.
These priors yield a particular collection
of posterior distributions $\Pi_\a(\cdot\mid X)$, and even the
best possible (or oracle) choice
of the tuning parameter $\a$ procedure is restricted to work through
this collection of
posterior distributions. Thus, the resulting procedure cannot be
expected to be optimal
for every model. One may think, for instance, of a model defined
through a wavelet expansion,
which has a double index and may not fit the Sobolev scale.
Second, even in a situation that the collection $\Pi_\a(\cdot\mid
X)$ contains
an optimal candidate, the empirical Bayes procedure (\ref{EqDefinitionHata}),
linked to the likelihood, although minimax over
the usual models, may fail to choose the optimal $\a$ for other models.
Other empirical Bayes procedures sometimes perform better, for
instance, by
directly relating to the bias--variance trade-off.

The radii $r_{n,\g}(\a)$ of the credible sets are decreasing in $\a$.
Hence, if the empirical Bayes choice $\hat\a_n$ in (\ref{EqDefinitionHata})
is restricted to a bounded interval
$[0,\Aupper]$, then the credible set $\hat C_n(L)$ has radius not
smaller than $r_{n,\g}(A)$,
which is bigger than necessary if the true parameter has greater
``regularity'' than $A$.
By the second statement of the theorem this can be remedied by choosing
$A$ dependent on $n$, at the cost of increasing the radius by a slowly
varying term.

\subsection{Hyperrectangles}\label{SectionHyperrectangle}
The \emph{hyperrectangle} $\Theta^\b(M)$
of order $\b$ and radius $M$ is defined in (\ref{defhyperrectangle}).
The minimax risk for this model in the case of the direct (not inverse) problem
where $\kappa_i=1$ is given in \citet{DonohoLiuMcGibbon}.
A slight variation of their proof gives that the minimax risk for
square loss in our problem\vspace*{1pt}
is bounded above and below by multiples of
$M^{{(1+2p)}/{(1+2\b+2p)}} n^{-{2\b}/({1+2\b+2p})}$, where the
constant depends on $C$
and $p$ in (\ref{assumpinverse}) only.
Furthermore, this order does not change if the hyperrectangle is reduced
to self-similar sequences.

%
\begin{proposition}\label{LemmaMinimaxRate}\label{LemmaComplexity}
Assume (\ref{assumpinverse}). For all $\b, M > 0$,
\[
\inf_{\hat{\th}_n}\sup_{\th_0\in\Theta^\b(M)}E_{\th_0}
\llVert\hat{\th}_n-\th_0\rrVert^2\asymp
M^{\vafrac{1+2p}{1+2\b+2p}} n^{-\afrac{2\b}{1+2\b+2p}},
\]
where the infimum is over all estimators.
This remains true if $\Theta^\b(M)$ is replaced
by $\Theta_{ss}^\b(M,\eps)$, for any sufficiently small $\eps>0$.
\end{proposition}

\begin{pf}
The problem of estimating $(\th_i)$ based on the data (\ref{eqdata})
is equivalent to estimating $(\th_i)$ based on independent $Y_1,
Y_2,\ldots$
with $Y_i\sim N(\th_i,\break n^{-1}\kappa_i^{-2})$.
As explained in \citet{DonohoLiuMcGibbon} (who consider identical variances
instead of $\s_i^2=n^{-1}\kappa_i^{-2}$ depending on $i$, but this
does not
affect the \mbox{argument}), the minimax estimator
for a given hyperrectangle is the vector of estimators
$T=(T_1,T_2,\ldots)$, where
$T_i$ is the minimax estimator in the problem of estimating $\th_i$
based on the single $Y_i\sim N(\th_i, \s_i^2)$, for each $i$, where
it is known that $\th_i^2\le M_i:=Mi^{-1-2\b}$. Furthermore,
\citet{DonohoLiuMcGibbon} show that in these univariate problems
the minimax risk when restricting to estimators $T_i(Y_i)$ that are
linear in $Y_i$
is at most $5/4$ times bigger than the (unrestricted, true) minimax
risk, where the former
\emph{linear minimax risk} is easily computed to be equal to $M_i\s
_i^2/(M_i+\s_i^2)$.
Thus, the minimax risk in the present
situation is up to a factor $5/4$ equal to
\[
\sum_{i=1}^\infty\frac{M_i\s_i^2}{M_i+\s_i^2} =\sum
_{i=1}^\infty\frac{i^{-1-2\b}Mn^{-1}\kappa_i^{-2}}{i^{-1-2\b
}M+n^{-1}\kappa_i^{-2}}.
\]
Using assumption (\ref{assumpinverse}) and Lemma~\ref{LemmaTechLem10}
(with $l=1$, $m=0$, $r=1+2\b+2p$, $s=2p$, and $Mn$ instead of $n$),
we can evaluate this as the right-hand side of the proposition.

To prove the final assertion, we note that
the self-similar functions $\Theta_{ss}^\b(M,\eps)$ are sandwiched between
$\Theta^\b(M)$ and the set
\[
\bigl\{\th\dvtx  \sqrt{\eps M\rho^{1+2\b}} i^{-1/2-\b} \leq\th
_i\leq\sqrt{M} i^{-1/2-\b}\mbox{ for all }i \bigr\}.
\]
Likewise, the minimax risk for $\Theta_{ss}^\b(M,\eps)$ is
sandwiched between the minimax risks of these models.
The smaller model, given in the display, is actually also a
hyperrectangle, which can be shifted to the
centered hyperrectangle $\Theta^{\b}((1-\sqrt{\eps\rho^{1+2\b
}})^2M/4)$. By shifting the observations
likewise we obtain an equivalent experiment. The $\ell^2$-loss is
equivariant under this shift.
Therefore, the minimax risks of the smaller and bigger models in the sandwich
are proportional to $(1-\sqrt{\eps\rho^{1+2\b}})^2
M^{({1+2p})/({1+2\b+2p})} n^{-2\b/(1+2\b+2p)}/4$
and $M^{({1+2p})/({1+2\b+2p})} n^{-2\b/(1+2\b+2p)}$, respectively,
by the preceding paragraph.
\end{pf}

Theorem~\ref{thmCoverage} shows that the credible sets $C_n(L)$
cover self-similar
parameters in $\Theta^\b(M)$
and, more generally, parameters satisfying the polished tail condition,
uniformly in regularity parameters $\b\in[0,\Aupper]$ and also uniformly
in the radius $M$ (but dependent on $\eps$ in the definition of
self-similarity or
$L_0$ in the definition of the polished tail condition).

Straightforward adaptation of the proof of Theorem~2 in \citet{KSzVZ}
shows that the empirical Bayes
posterior distribution $\Pi_{\hat\a_n}(\cdot\mid X)$ contracts to
the true parameter at the
minimax rate, by a logarithmic factor, uniformly over any
hyperrectangle $\Theta^\b(M)$ with
$\b\le\Aupper$. This immediately implies that the radius $r_{n,\g
}(\hat\a_n)$ is at most a
logarithmic factor larger than the minimax rate, and hence the size of
the credible sets $\hat
C_n(L)$ adapts to the scale of hyperrectangles. Closer inspection shows
that the logarithmic factor
in this result arises from the bias term and is unnecessary for the radius
$r_{n,\g}(\hat\a_n)$. Furthermore, this radius also adapts to the
constant $M$ in the optimal
manner.

%
\begin{proposition}\label{lemmaPropositionRectangleMinimxRateEB}
For every $\b\in(0,\Aupper]$ and $M > 0$,
\[
\inf_{\th_0 \in\Theta^\b(M)} P_{\th_0} \bigl( r_{n,\gamma}(\hat
\a_n)\le K M^{\vafrac{1/2+p}{1+2\b+2p}}n^{-\afrac{\b}{1+2\b+2p}} \bigr)
\rightarrow1.
\]
\end{proposition}

The proof of this proposition and of all further results in this section
is given in Section~\ref{SectionProofHyper}.

This encouraging result can be refined to an oracle type result
for self-similar true parameters. Recall that $\Theta_{ss}^\b(M,\eps)$
is the collection of self-similar parameters in $\Theta^\b(M)$.

%
\begin{theorem}
\label{TheoremRectangleRadiusSelfSimilar}
For all $\eps$ there exists a constant $K(\eps)$ such that, for all $M$,
%
%
\begin{equation}
\inf_{\th_0\in\Theta_{ss}^\b(M,\eps)}P_{\th_0} \Bigl(r_{n,\gamma
}^2(
\hat{\a}_n) \leq K(\eps) \inf_{\a\in[0,\Aupper]}
E_{\th_0}\llVert\hat\th_{n,\a
}-\th_0\rrVert
^2 \Bigr)\rightarrow1. \label{eqadaptivity}
\end{equation}
\end{theorem}

The theorem shows that for self-similar truths $\th_0$ the square
radius of the credible set
is bounded by a multiple of the mean square error\break
$\inf_{\a\in[0,\Aupper]} E_{\th_0}\llVert \hat\th_{n,\a}-\th_0\rrVert ^2$
of the best estimator in the class of
all Bayes estimators of the form $\hat\th_{n,\a}$, for $\a\in
[0,\Aupper]$.
The choice $\a$ can be regarded as made by an oracle
with knowledge of~$\th_0$. The class of estimators $\hat\th_{n,\a}$
is not complete, but
rich. In particular, the proposition below shows that it
contains a minimax estimator for every hyperrectangle $\Theta^\b(M)$
with $\b\le\Aupper$.

%
\begin{proposition}
\label{PropositionRectangleCompletenessBayesEstimators}
For every $\b\in(0,\Aupper]$ and $M > 0$,
\begin{eqnarray*}
\inf_{\a\in[0,\Aupper]} \sup_{\th_0\in\Theta^{\b}(M)}E_{\th
_0}
\llVert\hat\th_{n,\a}-\th_0\rrVert^2 &\lesssim& M^{\vafrac{1+2p}{1+2p+2\b}}n^{-\afrac{2\b}{1+2\b+2p}}.
\end{eqnarray*}
\end{proposition}

Within the scale of hyperrectangles it is also possible to study the
asymptotic behavior of the empirical Bayes regularity $\hat\a_n$
and corresponding posterior distribution.
For parameters in $\Theta_{ss}^\b(M)$ the empirical Bayes estimator estimates
$\beta$, which might thus be considered a ``true'' regularity $\b$.

%
\begin{lemma}\label{lemmaBoundsForAinHyperrectangle}
For any $0 < \b\le\Aupper-1$ and $M \ge1$, there exist constants
$K_1$ and $K_2$
such that $\Pr_{\th_0}(\b-K_1/\log n\le\hat\a_n\le\b+K_2/\log
n)\rightarrow1$ uniformly in
$\th_0\in\Theta_{ss}^\b(M,\eps)$.
\end{lemma}

Inspection of the proof of the lemma shows that the constant $K_2$ will
be negative
for large enough $M$, meaning that the empirical Bayes estimate $\hat
\a_n$ will
then actually slightly (up to a $1/\log n$ factor) undersmooth the
``true'' regularity.
The assumption that $\th_0\in\Theta_{ss}^\b(M,\eps)$ implies not
only that
$\th_0$ is of regularity $\b$, but also that within the rectangle it
is at a distance proportional
to $M$ from the origin. Thus, an increased distance from the origin
in a rectangle of fixed regularity $\b$ is viewed by the empirical
Bayes procedure as
``a little less smooth than $\b$.'' (This is intuitively reasonable
and perhaps the smoothness
of such $\th_0$ should indeed be viewed as smaller than $\beta$.)
The ``undersmoothing'' is clever, and actually
essential for the coverage of the empirical Bayes credible sets,
uniformly over all radii $M>0$.

\subsection{Sobolev balls}\label{SectionSobolevballs}
The \emph{Sobolev ball} $S^\b(M)$
of order $\b$ and radius $M$ is defined in (\ref{defsobolevspace}).
The minimax risk for this model is well known to be of the same
order as for the hyperrectangles considered in Section~\ref
{SectionHyperrectangle}
[see \citet{CavURE} or \citet{Cavalier2}]. A restriction to polished
tail sequences
decreases this by at most a logarithmic factor.

%
\begin{proposition}
\label{PropositionSobolevSize}
Assume (\ref{assumpinverse}). For all $\b, M > 0$,
\begin{eqnarray*}
&& \inf_{\hat{\th}_n}\sup_{\th_0\in S^\b(M)\cap\Theta
_{pt}(L_0)}E_{\th_0}
\llVert\hat{\th}_n-\th_0\rrVert^2
\\
&&\qquad \gtrsim
M^{\vafrac{1+2p}{1+2\b+2p}} n^{-\afrac{2\b}{1+2\b+2p}} (\log n)^{-\vafrac
{2+4p}{1+2\b+2p}},
\end{eqnarray*}
where the infimum is over all estimators.
\end{proposition}

\begin{pf}
The set $S^\b(M)\cap\Theta_{pt}(L_0)$ contains the set
\[
\bigl\{\th\in\ell_2\dvtx  \eps M i^{-1-2\b}(\log i)^{-2}\le\th_i\le\eps
' M i^{-1-2\b}(\log i)^{-2}\bigr\},
\]
for suitable $\eps<\eps'$. This is a translate of a hyperrectangle delimited
by the rate sequence $Mi^{-1-2\b}(\log i)^{-2}$. The minimax risk over this
set can be computed as in the proof of Proposition~\ref{LemmaMinimaxRate}.
\end{pf}

By Theorem~2 of \citet{KSzVZ} the empirical Bayes posterior distribution
$\Pi_{\hat\a_n}(\cdot\mid X)$ (with $\Aupper$ set to $\log n$) contracts
to the true parameter at the minimax rate, up to a logarithmic factor, uniformly
over any Sobolev ball $S^\b(M)$.
This implies that the radius $r_{n,\g}(\hat\a_n)$ is at most a
logarithmic factor
larger than the minimax rate, and hence the size of the credible sets
$\hat C_n(L)$ adapts
to the Sobolev scale. Again, closer inspection shows that
the logarithmic factors do not enter into the radii, and the radii
adapt to
$M$ in the correct manner. The proof of the following theorem can be found
in Section~\ref{SectionProofSobolev}.

%
\begin{theorem}
\label{TheoremSobolevRadiusPolishedTail}
There exists a constant $K$ such that, for all $\b\in[0,\Aupper]$,
$M$ and $L_0$,
\[
\inf_{\th_0\in S^\b(M)}P_{\th_0} \bigl(r_{n,\gamma}^2(
\hat{\a}_n) \leq K M^{\vafrac{1+2p}{1+2\b+2p}} n^{-\afrac{2\b}{1+2\b+2p}}
\bigr)
\rightarrow1.
\]
\end{theorem}

Thus, the empirical Bayes credible sets are honest confidence sets,
uniformly over the scale of Sobolev balls (for $\b\in[0,\Aupper]$
and $M>0$) intersected
with the polished tail sequences, of the minimax size over every
Sobolev ball.

If $\Aupper$ in (\ref{EqDefinitionHata}) is allowed to grow at the
order $\sqrt{\log n}$,
then these results are true up to logarithmic factors.

\subsection{Supersmooth parameters}\label{SectionSupersmooth}
Classes of \emph{supersmooth parameters} are defined by, for fixed $N,
M, c, d>0$,
%
%
\begin{eqnarray}
C^{00}(N_0,M)&=&\Bigl\{\th\in\ell^2\dvtx  \sup
_{i> N_0}\llvert\th_i\rrvert=0\mbox{ and }
\sup_{i\le N_0}\llvert\th_i\rrvert\leq
M^{1/2}\Bigr\},\label{defcnull}
\\
S^{\infty,c,d}(M)&=&\Biggl\{\th\in\ell^2\dvtx  \sum
_{i=1}^\infty e^{ci^d}\th_i^2
\le M\Biggr\}.\label{defsinfty}
\end{eqnarray}
The minimax rate over these classes is $n^{-1/2}$ or $n^{-1/2}$ up to a
logarithmic factor.
Every class $C_{00}(N_0,M)$ for fixed $N_0$ and arbitrary $M$ consists
of polished tail sequences only.

If the empirical Bayes regularity $\hat\a_n$ in (\ref
{EqDefinitionHata}) is restricted
to a compact interval $[0,\Aupper]$, then it will tend to the upper
limit $\Aupper$
whenever the true parameter is in one of these supersmooth models. Furthermore,
the bias of the posterior mean will be of order $O(1/n)$, which is negligible
relative to its variance at $\a=\Aupper$ (which is the smallest over
$[0,\Aupper]$).
Coverage by the credible sets then results.

More interestingly, for $\hat\a_n$ in (\ref{EqDefinitionHata}) restricted
to $[0,\log n]$, it is shown in \citet{KSzVZ} that the posterior
distribution
contracts at the rate $n^{-1/2}$ up to a lower order factor, for
$\th_0\in S^{\infty,c,1}$. Thus, the empirical Bayes credible sets then
adapt to the minimal size over these spaces. Coverage is not automatic,
but does take place uniformly in polished tail sequences by
Theorem~\ref{thmCoverage}.

The following theorem extends the findings on the sizes of the credible sets.

%
\begin{theorem}
\label{TheoremSupersmoothRadiusPolishedTail}
There exists a constant $K$ such that, for $\Aupper=\Aupper_n=\sqrt
{\log n}/(4\log\rho)$,
%
%
\begin{eqnarray}
\inf_{\th_0\in C^{00}(N_0,M)}P_{\th_0} \bigl(r_{n,\gamma}^2(
\hat{\a}_n) \leq K e^{(3/2+3p)\sqrt{\log N_0}\sqrt{\log n}} n^{-1}
\bigr)
\rightarrow1,\label{rateC00}
\\
\inf_{\th_0\in S^{\infty,c,d}(M)}P_{\th_0} \bigl(r_{n,\gamma
}^2(
\hat{\a}_n) \leq e^{(1/2+p)\sqrt{\log n}\log\log n} n^{-1} \bigr
)\rightarrow
1.\label{rateSinfty}
\end{eqnarray}
\end{theorem}

\section{Simulation example}
\label{secsim}
We investigate the uncertainty quantification of the empirical Bayes
credible sets in an example.
Assume that we observe the process
\begin{eqnarray*}
X_t&=&\int_0^t\int
_0^s \th_0(u) \,du \,ds+
\frac{1}{\sqrt{n}}B_t, \qquad t \in[0,1],
\end{eqnarray*}
where $B$ denotes a standard Brownian motion and $\th_0$ is the
unknown function of interest. It is
well known and easily derived that this problem is equivalent to~(\ref
{eqdata}), with $\th_{0,i}$
the Fourier coefficients of $\th_0$ relative to the eigenfunctions
$e_i(t) =
\sqrt{2}\cos(\pi(i-1/2)t)$ of the Volterra integral operator $K\th
(t)=\int_0^{t}\th(u)\,du$,
and $\k_i=1/(i-1/2)/\pi$ the
corresponding eigenvalues. In particular, $p=1$ in~(\ref
{assumpinverse}), that is, the
problem is mildly ill posed.

For various signal-to-noise ratios $n$ we simulate data from this model
corresponding to the true
function
%
%
\begin{equation}
\th_0(t)=\sum_{i=1}^{\infty}
\bigl(i^{-3/2}\sin(i)\bigr)\sqrt{2}\cos\bigl(\pi(i-1/2)t
\bigr).\label{exselfsim}
\end{equation}
This function $\th_0$ is self-similar with regularity parameter $\b=1$.
%
%
\begin{figure}

\includegraphics{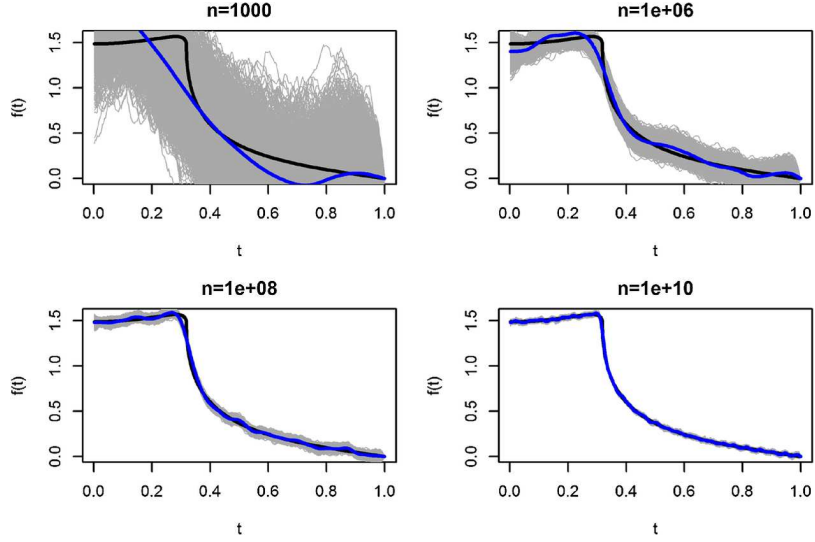}

\caption{Empirical Bayes credible sets. The true function is drawn in
black, the posterior mean in blue and the credible set in grey. We have
$n=10^3,10^6,10^8$ and $10^{10}$, respectively.}
\label{figure1}
\end{figure}
In Figure~\ref{figure1} we visualize $95\%$ credible sets for $\th_0$
(gray), the posterior mean
(blue) and the true function (black), by simulating $2000$ draws
from the empirical Bayes posterior distribution and plotting the $95\%$
draws out of the 2000
that are closest to
the posterior mean in the $L^2$-sense. The credible sets are drawn for
$n = 10^4, 10^6, 10^8$ and
$10^{10}$, respectively. The pictures show good coverage, as predicted
by Theorem~\ref{thmCoverage}. We note that we did not blow up the
credible sets by a
factor $L > 1$.

\begin{figure}

\includegraphics{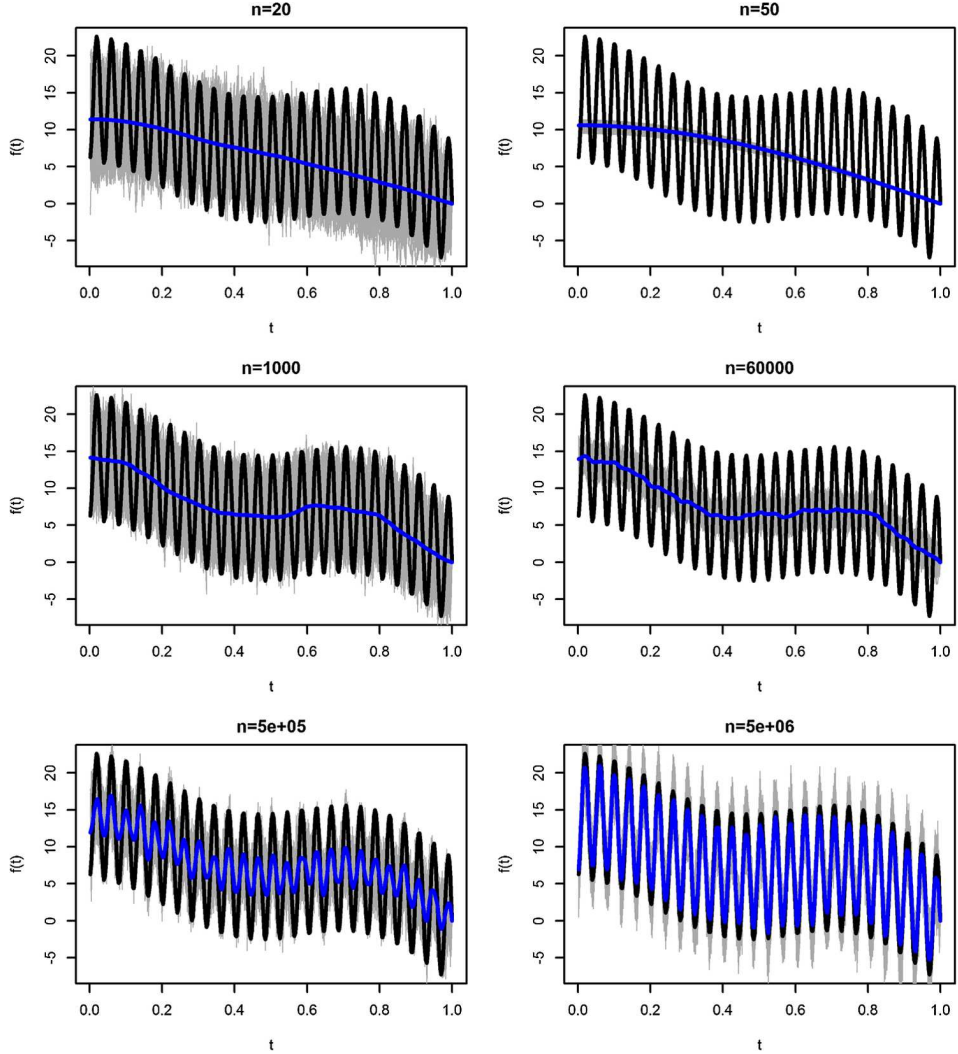}

\caption{Empirical Bayes credible sets for a non self-similar function.
The true function is drawn in black, the posterior mean in blue and the
credible sets in grey. From left to right and top to bottom we have
$n=20,50,10^3,6*10^4,5*10^5$ and $5*10^6$.}
\label{figureCounterExample}
\end{figure}

To illustrate the negative result of Theorem~\ref{LemmaCounterExample},
we also computed credible sets for a ``bad truth.'' We simulated data
using the following function:
\begin{eqnarray*}
&&\th_0(t)=\sum_{i=1}^{\infty}
\th_{0,i}\sqrt{2}\cos\bigl(\pi(i-1/2)t\bigr),
\end{eqnarray*}
where
\begin{eqnarray*}
\th_{0,i}=\cases{ 8, &\quad for $i=1$,
\vspace*{3pt}\cr
2, &\quad for $i=3$,
\vspace*{3pt}\cr
-2,
&\quad if $i=50$,
\vspace*{3pt}\cr
i^{-3/2}, &\quad if $2^{4^j}<i
\leq2*2^{4^j}$, for $j\geq3$,
\vspace*{3pt}\cr
0, &\quad else.}
\end{eqnarray*}
Figure~\ref{figureCounterExample} shows the results, with again the
true function $\th_0$
in black, the posterior mean in blue and the credible sets in gray.
The noise levels are $n=20,50,10^3,6*10^4,5*10^5$ and $5*10^6$, respectively.
As predicted by Theorem~\ref{LemmaCounterExample},
the coverage is very bad along a subsequence. For certain values of $n$
the posterior mean is far from the truth, yet the credible sets very
narrow, suggesting large
confidence in the estimator.

%

\section{Proof of Theorem~\texorpdfstring{\protect\ref{thmCoverage}}{3.6}}
\label{secproofcov}
\label{secMain}
The proof of Theorem~\ref{thmCoverage} is based on a
characterization of
two deterministic bounds on the data-driven choice $\hat\a_n$
of the smoothing parameter. Following \citet{KSzVZ} and \citet{Szabo},
define a function $h_n=h_n(\cdot; \th_0)$ by
%
%
\begin{equation}
h_n(\a;\th_0)= \frac{1+2\a+2p}{n^{1/(1+2\a+2p)}\log n} \sum
_{i=1}^{\infty}\frac{n^2i^{1+2\a}(\log i)\th_{0,i}^2}{
(i^{1+2\a+2p}+n)^2}\label{eqhnalpha}, \qquad\a
\ge0.
\end{equation}
Next define
%
%
\begin{eqnarray}
\underline{\a}_n(\th_0) &=&\inf\bigl\{\a\in[0,
\Aupper]\dvtx  h_n(\a;\th_0)\geq1/\bigl(16C^8
\bigr)\bigr\},\label
{defLowAlpha}
\\
\overline{\a}_n(\th_0) &=&\sup\bigl\{\a\in[0,
\Aupper]\dvtx  h_n(\a;\th_0)\leq8C^8 \bigr
\}\label{defUppAlpha}.
\end{eqnarray}
An infimum or supremum over an empty set can be understood to be
$\Aupper$ or 0, respectively.
The value $\Aupper$ is as in (\ref{EqDefinitionHata}). If it depends
on $n$,
then this dependence is copied into the definitions.

The following theorem shows that uniformly in parameters $\th_0$ that
satisfy the polished
tail condition, the sequences $\underline\a_n(\th_0)$ and $\overline
\a_n(\th_0)$ capture $\hat\a_n$ with probability
tending to one. Furthermore, these sequences are at most a universal
multiple of $1/\log n$
(or\vspace*{1pt} slightly more if $\Aupper$ depends on $n$) apart, leading to the same
``rate'' $n^{-\a/(1+2\a+2p)}$, again uniformly in polished tail
sequences $\th_0$.

%
\begin{theorem}\label{thmmle}
For every $L_0\ge1$,
%
%
\begin{eqnarray}
\inf_{\th_0\in\Theta_{pt}(L_0)} P_{\th_0} \bigl(\underline
\a_n(\th_0)\le\hat\a_n\le\overline\a
_n(\th_0) \bigr)&\rightarrow& 1,\label{eqm1}
\\
\sup_{\th_0\in\Theta_{pt}(L_0)}\frac{n^{-2\underline\a_n(\th
_0)/(1+2\underline\a_n(\th_0)+2p)}}{
n^{-2\overline\a_n(\th_0)/(1+2\overline\a_n(\th_0)+2p)}}&\le&
K_{3,n},\label{eqDiffBounds}
\end{eqnarray}
for $\log n\ge2+4\Aupper+2p\vee C_0$, where $C_0$ depends on $N_0$,
\begin{eqnarray*}
&& K_{3,n}\le c \bigl(2^9C^{16}L_0^2
\rho^{5+10\Aupper+6p}\bigr) ^{1+2p},
\end{eqnarray*}
for some universal constant $c$.
\end{theorem}

\begin{pf}
The proof of (\ref{eqm1}) is similar to proofs given in \citet{KSzVZ}.
However, because the exploitation of the polished tail condition and
the required uniformity in the parameters are new, we provide
a complete proof of this assertion in the \hyperref[secproofMLE2]{Appendix}.
Here we only prove inequality (\ref{eqDiffBounds}).

Let $N_\a=n^{1/(1+2\a+2p)}$ and set
\[
\Theta_j=\sum_{i=\rho^{j-1}+1}^{\rho^j}
\th_i^2.
\]
Then the polished tail condition (\ref{EqPolishedTail}) implies that
$\Theta_j\le L_0 \Theta_{j'}$ for all $j\ge j'$, whence
\begin{eqnarray*}
h_n(\a;\th) &=& \frac{1}{N_\a\log N_\a}\sum
_{i=1}^\infty\frac{n^2i^{1+2\a}\th
_i^2\log i}{(i^{1+2\a+2p}+n)^2}
\\
&\le&\frac{1}{N_\a\log N_\a}\sum_{j=1}^\infty
\Theta_j\frac
{n^2\rho^{j(1+2\a)}j\log\rho}{
(\rho^{(j-1)(1+2\a+2p)}+n)^2}.
\end{eqnarray*}
We get a lower bound if we swap the $j$ and $j-1$ between numerator and
denominator. [The term $(j-1)\log\rho$ that then results in the
numerator is
a nuisance for $j=1$; instead of $(j-1)\log\rho= \log\rho^{j-1}$,
we may use $\log2$
if $j=1$, as the term $i=1$ does not contribute.]

Define $J_\a$ to be an integer such that
\begin{eqnarray}
J_\a= \biggl\lfloor\frac{\log n}{(\log\rho)(1+2\a+2p)} \biggr\rfloor\nonumber
\\
\eqntext{\mbox{hence } \rho^{J_\a(1+2\a+2p)}\leq n< \rho^{1+2\a
+2p}\rho^{J_\a(1+2\a+2p)}.}
\end{eqnarray}

Under the polished tail condition,
\begin{eqnarray*}
&& \sum_{j>J_\a}\Theta_j\frac{n^2\rho^{j(1+2\a)}j\log\rho} {(\rho
^{(j-1)(1+2\a+2p)}+n)^2}
\\
&&\qquad \le  L_0\Theta_{J_\a}\sum_{j>J_\a}n^2
\rho^{-j(1+2\a+4p)}j(\log\rho)\rho^{2(1+2\a+2p)}
\\
&&\qquad \lesssim L_0\Theta_{J_\a}n^2
\rho^{-J_\a(1+2\a+4p)}J_\a(\log\rho) \rho^{2(1+2\a+2p)}
\\
&&\qquad \leq L_0\Theta_{J_\a}nJ_\a(\log\rho)
\rho^{3(1+2\a+2p)-J_{\a}2p}.
\end{eqnarray*}
The constant in $\lesssim$ is universal.
[We have
$\sum_{j>J} jx^{j}=x^{J+1}(J+1) ((1-x)^{-1}-(1-x)^{-2}/(J+1)
)$ for $0<x<1$,
and $x=\rho^{-(1+2\a+4p)}\le\rho^{-1}\le1/2$.]

The $J_\a$th term of the series on the left-hand side is
\[
\Theta_{J_\a}\frac{n^2\rho^{J_\a(1+2\a)}J_\a\log\rho} {(\rho
^{(J_\a-1)(1+2\a+2p)}+n)^2} \geq(1/4) \Theta_{J_\a} n
J_\a\rho^{-J_{\a}2p}\log\rho.
\]
Up to a factor $L_0\rho^{3(1+2\a+2p)}$ this
has the same order of magnitude as the right-hand side of the preceding
display, whence
\[
h_n(\a;\th) \lesssim\bigl(1+L_0\rho^{3+6\a+6p}
\bigr) \frac{1}{N_\a\log N_\a} \sum_{j=1}^{J_\a}
\Theta_j\frac{n^2\rho^{j(1+2\a)}j(\log\rho)}{
(\rho^{(j-1)(1+2\a+2p)}+n)^2}.
\]

Because $\rho^{j(1+2\a+2p)}+n\le2n$ for $j\le J_\a$, we also have
\begin{eqnarray*}
h_n(\a;\th) &\ge&\frac{1}{N_\a\log N_\a}\sum
_{j=1}^{J_\a}\frac{\Theta_j\rho
^{(j-1)(1+2\a)}\log(\rho^{j-1}\vee2)}{4}
\nonumber
\\
&\ge&\frac{1}{N_\a\log N_\a}\frac{\rho^{-(1+2\a)}}8\sum_{j=1}^{J_\a}
\Theta_j\rho^{j(1+2\a)}j(\log\rho\wedge1).
\end{eqnarray*}
(Note that $\log2\approx0.69\ge1/2$ and $j-1\ge j/2$ for $j\ge2$.)\vspace*{1pt}

Now fix $\a_1\le\a_2$. Then $J_{\a_2}\le J_{\a_1}$ and $\rho
^{1+2\a_1}\le\rho^{1+2\a_2}$ and
\begin{eqnarray*}
&& \frac{h_n(\a_1;\th)}{h_n(\a_2;\th)}
\\
&&\qquad \lesssim \frac{N_{\a_2}\log N_{\a_2}}{N_{\a_1}\log N_{\a_1}} \frac
{(1+L_0\rho^{3+6\a_1+6p}) (\sum_{j=1}^{J_{\a_2}}+\sum_{j=J_{\a
_2}}^{J_{\a_1}})\Theta_j\rho^{j(1+2\a_1)}j(\log\rho)}{
\rho^{-(1+2\a_2)}\sum_{j=1}^{J_{\a_2}}\Theta_j\rho^{j(1+2\a
_2)}j(\log\rho\wedge1)}
\\
&&\qquad \lesssim \frac{N_{\a_2}\log N_{\a_2}}{N_{\a_1}\log N_{\a
_1}}\bigl(1+L_0\rho^{3+6\a_1+6p}\bigr)
\rho^{1+2\a_2} \biggl(1+\frac{L_0\Theta_{J_{\a_2}}\sum_{j=J_{\a
_2}}^{J_{\a_1}}\rho^{j(1+2\a_1)}j}{
\Theta_{J_{\a_2}}\rho^{J_{\a_2}(1+2\a_2)}J_{\a_2}} \biggr)
\\
&&\qquad \lesssim\frac{N_{\a_2}\log N_{\a_2}}{N_{\a_1}\log N_{\a
_1}}\bigl(1+L_0\rho^{3+6\a_1+6p}\bigr)
\rho^{1+2\a_2} \biggl(1+L_0\frac{\rho^{J_{\a_1}(1+2\a_1)}J_{\a_1}}{
\rho^{J_{\a_2}(1+2\a_2)}J_{\a_2}} \biggr)
\\
&&\qquad \leq\frac{N_{\a_2}\log N_{\a_2}}{N_{\a_1}\log N_{\a
_1}}\bigl(1+L_0\rho^{3+6\a_1+6p}\bigr)
\rho^{2+4\a_2} \biggl(1+L_0\frac{J_{\a_1}}{J_{\a_2}} \biggr)
\\
&&\qquad \leq\frac{N_{\a_2}}{N_{\a_1}}\frac{1+2\a_1+2p}{1+2\a
_2+2p}\bigl(1+L_0
\rho^{3+6\a_1+6p}\bigr) \rho^{2+4\a_2} \biggl(1+2L_0
\frac{1+2\a_2+2p}{1+2\a_1+2p} \biggr).
\end{eqnarray*}
[We use $\sum_{j=I}^J x^j=x^I(x^{J-I+1}-1)/(x-1)\le x^J
x/(x-1)\lesssim x^J$, for $x=\rho^{1+2\a}\ge\rho$.]

Since $\underline\a_n\le\overline\a_n$,
there is nothing to prove in the trivial cases
$\underline\a_n=\Aupper$ or $\overline\a_n=0$.
In the other cases it follows that $h_n(\underline\a_n;\th_0)\ge
1/(16C^8)$ and $h_n(\overline\a_n;\th_0)\le8 C^8$.
Then the left-hand side of the preceding display with $\a_1=\underline\a
_n$ and $\a_2=\overline\a_n$
is bounded from below by $1/(128C^{16})$. After taking the $(1+2p)$th power
of both sides and\vadjust{\goodbreak} rearranging the inequality we get
\begin{eqnarray*}
&& \bigl(2^8C^{16}\bigr)^{1+2p}\bigl(1+L_0
\rho^{3+6\Aupper+6p}\bigr)^{1+2p}\rho^{(2+4\Aupper)(1+2p)}(1+2L_0)^{1+2p}
\\
&&\qquad \gtrsim N_{\underline\a_n}^{1+2p}/N_{\overline\a_n}^{1+2p}
=N_{\underline\a_n}^{-2\underline\a_n}/N_{\overline\a
_n}^{-2\overline\a_n}.
\end{eqnarray*}
This concludes the proof of Theorem~\ref{thmmle}.
\end{pf}

We proceed to the proof of Theorem~\ref{thmCoverage}. Recall the
definition of the posterior distribution $\Pi_\a(\cdot\mid X)$ in
(\ref{PostDist}).

For notational convenience denote the mean of the posterior distribution
(\ref{PostDist}) by $\hat\th_\a$ and
the radius $r_{n,\gamma}(\a)$ defined by (\ref{eqradius}) by $r(\a
)$. Furthermore, let
%
%
\begin{equation}
\label{EqDefWB} W(\a) = \hat\th_\a- E_{\th_0}\hat
\th_\a\quad\mbox{and}\quad B(\a;\th_0) =
E_{\th_0}\hat\th_\a- \th_0
\end{equation}
be the centered posterior mean and the bias of the posterior mean, respectively.
The radius $r(\a)$ and the distribution of the variable $W(\a)$ under
$\th_0$
are free of $\th_0$. On the other hand, the bounds $\underline\a_n$
and $\overline\a_n$
do depend on $\th_0$, but we shall omit this from the notation.
Because the radius of
the credible set is defined to be infinite if $\hat\a_n=0$, it is not
a loss of generality
to assume that $\overline\a_n>0$. For simplicity we take $\Aupper$ in
(\ref{EqDefinitionHata}) independent of $n$,
but make the dependence of constants on $\Aupper$ explicit, so that
the general
case follows by inspection.

We prove below that there exist positive parameters $ C_1, C_2, C_3$ that
depend on $C, \Aupper, p, L_0, \rho$ only such that, for all $\th_0
\in\Theta_{pt}(L_0)$,
%
%
\begin{eqnarray}
\label{r} \inf_{\underline\a_n \le\a\le\overline\a_n}r(\a) &\ge& C_1
n^{-\afrac{\overline\a_n}{1+2\overline\a_n + 2p}},
\\
\label{b} \sup_{\underline\a_n \le\a\le\overline\a_n}\bigl\llVert B(\a;
\th_0)\bigr\rrVert&\le& C_2 n^{-\afrac{\overline\a_n}{1+2\overline\a_n + 2p}},
\\
\label{w} \inf_{\th_0 \in\Theta_{pt}(L_0)}P_{\th_0} \Bigl(\sup
_{\underline
\a_n \le\a\le\overline\a_n} \bigl\llVert W(\a)\bigr\rrVert&\le&
C_3n^{-\afrac{\overline\a_n}{1+2\overline\a_n + 2p}}
\Bigr) \rightarrow1.
\end{eqnarray}
We have $\th_0\in\hat C_n(L)$ if and only if $\llVert \hat\th_{\hat\a
}-\th_0\rrVert \le L r(\hat\a)$,
which is implied by $\llVert W(\hat\a)\rrVert \le Lr(\hat\a)-\llVert
B(\hat\a;\th
_0)\rrVert $, by the
triangle inequality. Consequently, by (\ref{eqm1}) of Theorem~\ref
{thmmle},
to prove (\ref{eqPostCov}), it suffices to show that for $L$ large enough
\[
\inf_{\th_0 \in\Theta_{pt}(L_0)}P_{\th_0} \Bigl(\sup_{\underline
\a_n \le\a\le\overline\a_n}
\bigl\llVert W(\a)\bigr\rrVert\le L \inf_{\underline\a_n \le\a\le
\overline\a_n}r(\a) - \sup
_{\underline\a_n \le\a\le\overline\a_n}\bigl\llVert B(\a;\th_0)\bigr
\rrVert\Bigr)
\rightarrow1.
\]
This results from the combination of (\ref{r}), (\ref{b}) and (\ref
{w}), for $L$ such
that $C_3\le LC_1-C_2$.

We are left to prove (\ref{r}), (\ref{b}) and (\ref{w}).

\begin{pf*}{Proof of Proof of (\ref{r})}
The radius $r(\a)$ is determined by the requirement that
$P (U_n(\a)<r^2(\a) )=1-\gamma$,
for the random variable $U_n(\a)=\sum_i s_{i,n,\a}Z_i^2$,
where $s_{i,n,\a}=\k_i^{-2}/(i^{1+2\a}\k_i^{-2}+n)$
and $(Z_i)$ is an i.i.d. standard normal sequence. Because $\a\mapsto
s_{i,n,\a}$
is decreasing in $\a$, the map $\a\mapsto r(\a)$ is nonincreasing, and
hence the infimum in (\ref{r}) is equal to $r(\overline\a_n)$.
In view of (\ref{assumpinverse}) and Lemma~\ref{LemmaTechLem10},
the expected value and variance of $U_n(\a)$ satisfy, for $n\ge
e^{1+2\a+2p}$,
\begin{eqnarray*}
\E U_n(\a)&=&\sum_{i=1}^{\infty}s_{i,n,\a}
\geq\frac{1}{C^{4}}\sum_{i=1}^{\infty}
\frac{i^{2p}}{i^{1+2\a+2p}+n}
\\
&\ge&\frac{1}{C^4(3^{1+2\a+2p}+1)} n^{-\afrac{2\a}{1+2\a+2p}},
\\
\Var U_n(\a)&=&2\sum_{i=1}^{\infty}s_{i,n,\a}^2
\leq2C^{8}\sum_{i=1}^{\infty}
\frac{i^{4p}}{(i^{1+2\a+2p}+n)^2}
\\
&\le& 10 C^8n^{-\vafrac{1+4\a}{1+2\a+2p}}.
\end{eqnarray*}
We see that the standard deviation of $U_n(\a)$ is negligible compared
to its expected value. This implies that all quantiles of $U_n(\a)$ are
of the order $\E U_n(\a)$. More precisely, by Chebyshev's inequality
we have that $\Pr(U<r^2)=1-\gamma$ implies that
$\E U-(1-\gamma)^{-1/2}\sd U\le r^2\le\E U+\gamma^{-1/2} \sd U$.
For $\a\le\Aupper$ the expectation $\E U_n(\a)$ is further bounded
below by $C_{1,1} n^{-2\a/(1+2\a+2p)}$, for
$C_{1,1}=C^{-4}(3^{1+2\Aupper+2p}+1)^{-1}$. Furthermore,
$\sd U_n(\a)$ is bounded above by $C_{1,2,n}n^{-2\a/(1+2\a+2p)}$, for
$C_{1,2,n}=\sqrt{10}C^4 n^{-(1/2)/(1+2\Aupper+2p)}\rightarrow0$.
Hence, for large enough $n$ we have $C_{1,1}/2\geq(1-\gamma
)^{-1/2}C_{1,2,n}$, whence
(\ref{r}) holds for $C_1^2= C_{1,1}/2$ and sufficiently large $n$.
If $\Aupper$ depends on $n$, then so does $C_{1,1}=C_{1,1,n}$, but
since $\Aupper\leq\sqrt{\log n}/4$ by assumption, we still
have that $C_{1,2,n} \ll C_{1,1,n}$, and the preceding argument
continues to apply.
\end{pf*}

\begin{pf*}{Proof of (\ref{b})} In view of the explicit expression for
$\hat\th_\a$ and (\ref{assumpinverse}),
%
%
\begin{equation}
\label{EqBiasSquare} \bigl\llVert B(\a;\th_0)\bigr\rrVert^2
=\sum_{i=1}^{\infty}\frac{i^{2+4\a}\k_i^{-4}\th
^2_{0,i}}{(i^{1+2\a}\k_i^{-2}+n)^2} \le
C^8\sum_{i=1}^{\infty}
\frac{i^{2+4\a+4p}\th^2_{0,i}}{(i^{1+2\a
+2p}+n)^2}.
\end{equation}
The first term of the sum is zero following from our assumption $\theta
_{0,1}=0$. Since the right-hand side of the preceding display is
(termwise) increasing in $\a$,
the supremum over $\a$ is taken at $\overline\a_n$.
Because the map $x\mapsto x^{-1}\log x$ is decreasing for $x\ge e$ and
takes equal values at $x=2$ and $x=4$, its minimum over $[2,N^{1+2\a
+4p}]$ is taken
at $N^{1+2\a+4p}$ if $N^{1+2\a+4p}\ge4$. Therefore, for $2\le i\le N$
we have that $i^{1+2\a+4p}\le N^{(1+2\a+4p)}\log i/\log N$ if $N\ge
4^{1/(1+2\a+4p)}$.
Applied with $N_\a=n^{1/(1+2\a+2p)}$, this shows that, for $n\ge
4^{(1+2\a+2p)/(1+2\a+4p)}$,
\[
\sum_{2\le i \le N_\a}\frac{i^{2+4\a+4p}\th^2_{0,i}}{(i^{1+2\a
+2p}+n)^2} \le
n^{-{2\a}/({1+2\a+ 2p})}h_n(\a;\th_0).
\]
This bounds the initial part of the series in the right-hand side of (\ref
{EqBiasSquare}).
For $\th_0\in\Theta_{pt}(L_0)$, the remaining part can be bounded
above by
\begin{eqnarray*}
&&\sum_{i\geq N_\a}\theta_{0,i}^2
\le L_{0}\sum_{i=N_\a}^{\rho N_\a}
\theta_{0,i}^2 
\leq L_{0}\bigl(
\rho^{1+2\a+2p}+1\bigr)^2 h_n(\a;
\th_0)N_\a^{-2\a},
\end{eqnarray*}
as is seen by lower bounding the series $h_n(\a;\th_0)$ by the sum of
its terms
from $N_\a$ to $\rho N_\a$.
Using the inequality $h_n(\overline\a_n;\th_0)\leq8C^{8}$, we can
conclude that
\[
\bigl\llVert B_n(\overline\a_n;\theta_0)
\bigr\rrVert^2 \leq\bigl(L_02\rho^{2+4\Aupper+4p}+1
\bigr) 8C^{8} n^{-2\overline\a
_n/(1+2\overline\a_n+2p)}.
\]
This concludes the proof of (\ref{b}), with $C_2^2=(L_02\rho
^{2+4\Aupper+4p}+1) 8C^{8}$.
\end{pf*}

\begin{pf*}{Proof of (\ref{w})}
Under $P_{\th_0}$ the variable $V_n(\a) = \llVert W(\a)\rrVert ^2$ is
distributed as $\sum t_{i,n,\a} Z^2_i$,
where $t_{i,n,\a}=n\k_i^{-2}/(i^{1+2\a}\k_i^{-2}+n)^2$
and $Z_i:=\sqrt n(X_i-E_{\th_0} X_i)$ are independent standard normal
variables.
This representation gives a coupling over $\a$ and, hence, $\sup
_{\underline\a_n\le\a\le\overline\a_n}\llVert W(\a)\rrVert ^2$
is distributed as\break $\sup_{\underline\a_n\le\a\le\overline\a
_n}\sum t_{i,n,\a} Z^2_i
=\sum t_{i,n,\underline\a_n} Z^2_i$, since the coefficients
$t_{i,n,\a}$ are decreasing in $\a$.
By Lemma~\ref{LemmaTechLem10},
%
%
\begin{eqnarray}\label{eqVn}
\E V_n (\a) &=&\sum_{i=1}^{\infty}t_{i,n,\a}
\le C^6\sum_{i=1}^\infty
\frac{ni^{2p}}{(i^{1+2\a+2p}+n)^2} \nonumber
\\
&\le&  C^65 n^{-\afrac{2\a}{1+2\a +2p}},
\nonumber\\[-8pt]\\[-8pt]\nonumber
\Var V_n(\a)&=&2\sum_{i=1}^{\infty}t_{i,n,\a}^2
\le2C^{12}\sum_{i=1}^{\infty}
\frac{n^2i^{4p}}{(i^{1+2\a+2p}+n)^4}
\\
& \le& 10 C^{12}n^{-\vafrac{1+4\a}{1+2\a+2p}}.
\nonumber
\end{eqnarray}
Again, the standard deviation of $V_n(\a)$ is of smaller order than
the mean.
By reasoning as for the proof of (\ref{r}), we
obtain (\ref{w}) with the constant $\sqrt{6C^6}$, but with the rate
$n^{-\a/(1+2\a+2p)}$
evaluated at $\underline\a_n$ rather than $\overline\a_n$. Although
the last one is smaller,
it follows by (\ref{eqDiffBounds}) that the square rates are
equivalent up to
multiplication by $K_{3,n}$ (which is fixed if $\Aupper$ is fixed).
Thus, (\ref{w}) holds with $C_3^2=6C^6K_{3,n}$.
\end{pf*}

%
\begin{remark}
In the proof of (\ref{b}) we used the assumption that the first
coordinate $\th_{0,1}$
of the true parameter is zero. As $h_n(\a;\th_0)$ does not depend on this
coefficient, it would otherwise be impossible to use this function
in a bound on the square bias and thus relate the bias to $\hat\a_n$.
\end{remark}

\section{Proofs for hyperrectangles}
\label{SectionProofHyper}
In this section we collect proofs for
the results in Section~\ref{SectionHyperrectangle}.
Throughout the section we set $N_\a=n^{1/(1+2\a+2p)}$
and use the abbreviations of the preceding section.

\subsection{Proof of Theorem \texorpdfstring{\protect\ref{TheoremRectangleRadiusSelfSimilar}}{3.9}}

By (\ref{eqm1}) the infimum in (\ref{eqadaptivity}) is bounded
from below by
\[
\inf_{\th_0\in\Theta^{\b}_{ss}(M)}P_{\th_0} \Bigl(\sup_{\underline\a_n
\le\a\le\overline\a_n}
r(\a) \leq K \inf_{\a\in[0,\Aupper]} E_{\th_0}\llVert\hat
\th_{\a}-\th_0\rrVert^2 \Bigr) - o(1).
\]
Here the probability is superfluous because $r(\a)$, $\underline\a
_n$ and $\overline\a_n$ are
deterministic. We separately bound the supremum and infimum inside the
probability (above and
below) by a multiple of $M^{(1+2p)/(1+2\beta+2p)}n^{-2\beta/(1+2\beta+2p)}$.

As was argued in the proof of (\ref{r}), the radius $r(\a)$ is nonincreasing
in $\a$ and, hence, its supremum is $r(\underline\a_n)$. Also,
similarly as in the proof of (\ref{r}),
but now using the upper bound
\[
\E U_n(\a) \le{C^{4}}\sum
_{i=1}^{\infty}\frac{i^{2p}}{i^{1+2\a+2p}+n} \le
C^4(3+2/\a) n^{-\afrac{2\a}{1+2\a+2p}},
\]
by Lemma~\ref{LemmaTechLem10}, we find that
%
%
\begin{equation}
\sup_{\a\in[\underline\a_n,\overline\a_n]}r(\a)^2\le C^4(3+2/
\underline\a_n) n^{-2\underline\a_n/(1+2\underline\a
_n+2p)}.\label{eqUBr}
\end{equation}
The sequence $\underline\a_n$ tends to $\b$, uniformly in $\th_0\in
\Theta^\b(M)$
by (\ref{eqBoundsForBounds1}) (below), whence $(3+2/\underline\a
_n)\leq(4+2/\beta)$.
A second application of (\ref{eqBoundsForBounds1}), also involving
the precise definition of the constant $K_1$, shows that\vspace*{1pt}
the preceding display is bounded above
by a multiple of $M^{(1+2p)/(1+2\beta+2p)}n^{-2\beta/(1+2\beta+2p)}$.

With the notation as in (\ref{EqDefWB}), the minimal mean square error
in the right-hand side of the probability is equal to
%
%
\begin{equation}
\inf_{\a\in[0,\Aupper]} E_{\th_0}\llVert\hat\th_{n,\a}-
\th_0\rrVert^2 =\inf_{\a\in[0,\Aupper]} \bigl[
\bigl\llVert B(\a;\theta_0)\bigr\rrVert^2+E_{\th_0}
\bigl\llVert W(\a)\bigr\rrVert^2 \bigr]. \label{eqOracle}
\end{equation}
The square bias and variance terms in this expression are given
in (\ref{EqBiasSquare}) and (\ref{eqVn}).
By (\ref{assumpinverse}) the square bias is bounded below by
%
%
\begin{equation}
\frac{1}{C^8}\sum_{i=1}^{\infty}
\frac{i^{2+4\a+4p}\th
^2_{0,i}}{(i^{1+2\a+2p}+n)^2} \geq\frac{1}{4C^8} \sum_{i=N_\a}^{\rho
N_\a}
\theta_{0,i}^2 \geq\frac{\eps M}{4C^{8}} N_{\a}^{-2\beta},\label{eqOracleBias}
\end{equation}
for $\th_0\in\Theta^{\b}_{ss}(M)$.
By Lemma \ref{LemmaTechLem10} the variance
term $E_{\th_0}\llVert W(\a)\rrVert ^2$ in (\ref{eqOracle}) is
bounded from below by
%
%
\begin{equation}
\frac{1}{C^6}\sum_{i=1}^\infty
\frac{ni^{2p}}{(i^{1+2\a+2p}+n)^2} \geq\frac{1}{C^6}\bigl(3^{1+2\Aupper+2p}+1
\bigr)^{-2} N_{\a}^{-2\a}.\label
{eqOracleVariance}
\end{equation}
The square bias is increasing in $\a$, whereas the variance is decreasing;
the same is true for their lower bounds as given.
It follows that their sum is bounded below by the height at which the
two curves
intersect. This intersection occurs for $\a$ solving
%
%
\begin{equation}
\frac{\eps M}{4C^2}N_{\a}^{-2\beta}=\bigl(3^{1+2\Aupper+2p}+1
\bigr)^{-2} N_{\a}^{-2\a}. \label{defao}
\end{equation}
Write the solution as $\a=\beta-K/\log n$, in terms of some parameter $K$
(which may depend on $n$ as does the solution). By elementary algebra,
we get that
%
%
\begin{eqnarray}\label{eqratesrewrite}
\qquad N_\a^{-2\a}&\equiv& n^{-\afrac{2\a}{1+2\a+2p}}\nonumber
\\
&=& n^{-\afrac{2\beta}{1+2\beta+2p}}
\bigl(e^{\afrac{2K}{1+2\beta-2K/\log n+2p}} \bigr)^{\vafrac{1+2p}{1+2\beta
+2p}},
\nonumber\\[-8pt]\\[-8pt]\nonumber
N_\a^{-2\b}&\equiv& n^{-\afrac{2\beta}{1+2\a+2p}}
\\
&=& n^{-\afrac{2\beta}{1+2\beta+2p}}
\bigl(e^{\afrac{2K}{1+2\beta-2K/\log n+2p}} \bigr)^{-\afrac{2\beta
}{1+2\beta+2p}}.
\nonumber
\end{eqnarray}
Dividing the first by the second, we see from (\ref{defao}) that
$K$ is the solution to
%
%
\begin{equation}
e^{\afrac{2K}{1+2\beta-2K/\log n+2p}}=\bigl(3^{1+2A+2p}+1\bigr)^2
\frac{\eps
M}{4C^2}.\label{eqdefK}
\end{equation}
Furthermore,\vspace*{1pt} (\ref{eqratesrewrite}) shows that the value
of the right-hand side of (\ref{eqOracle}) at the corresponding
$\a=\b-K/\log n$ is equal\vspace*{1pt} to a constant times\break
$M^{(1+2p)/(1+2\b
+2p)}  n^{-2\b/(1+2\b+2p)}$,
where the constant multiplier depends only on $\beta,\eps$ and
$\Aupper$.

\subsection{Proof of Proposition \texorpdfstring{\protect\ref{PropositionRectangleCompletenessBayesEstimators}}{3.10}}
In view of (\ref{eqOracle}), (\ref{EqBiasSquare}) and (\ref{eqVn})
and (\ref{assumpinverse}),
\[
\inf_{\a\in[0,\Aupper]} E_{\th_0}\llVert\hat\th_{n,\a}-
\th_0\rrVert^2\leq C^8\sum
_{i=1}^{\infty}\frac{i^{2+4\a+4p}\th
^2_{0,i}}{(i^{1+2\a+2p}+n)^2} + C^{6}\sum
_{i=1}^\infty\frac{ni^{2p}}{(i^{1+2\a+2p}+n)^2}.
\]
By Lemma \ref{LemmaTechLem10} and the definition of the
hyperrectangle $\Theta^{\beta}(M)$, one can see that the right-hand
side of the preceding display for $\a\leq\beta$ is bounded from
above by
%
%
\begin{equation}
5C^{8}MN_\a^{-2\beta}+5C^6N_\a^{-2\a}.\label{eqOracleRate}
\end{equation}
Then choosing $\a=\beta-K/\log n$ with constant $K$ satisfying
\[
e^{\afrac{2K}{1+2\beta-2K/\log n+2p}}=M,
\]
we get from (\ref{eqratesrewrite}) that (\ref{eqOracleRate}) is
bounded above by the constant times
\[
M^{(1+2p)/(1+2\b+2p)}n^{-2\b
/(1+2\b+2p)}.
\]

\subsection{Proof of Lemma \texorpdfstring{\protect\ref{lemmaBoundsForAinHyperrectangle}}{3.11}}
We show that
%
%
\begin{eqnarray}
\inf_{\theta_0\in\Theta^{\beta}(M)}\underline\a_n(\th_0)&
\geq& \b- K_1/\log n\label{eqBoundsForBounds1},
\\
\sup_{\theta_0\in\Theta^{\beta}_{ss}(M,\eps)}\overline\a_n(\th_0)&
\leq&\b+K_{2}/\log n\label{eqBoundsForBounds2}
\end{eqnarray}
hold for $\log n\ge(4K_1)\vee(\log N_0)^2\vee4(1+2p)$ and constants
$K_1, K_{2}$ satisfying
%
%
\begin{eqnarray}
\label{EqKs} e^{2K_1/(1+2\b-2K_1/\log n+2p)}&=&288Me^{2\Aupper+3}C^{8},
\\
\eps Me^{2K_{2}/(1+2\Aupper+2p)}&=&\bigl(\rho^{1+2\Aupper+2p}+1\bigr)^28C^8.
\end{eqnarray}

\begin{pf*}{Proof of (\ref{eqBoundsForBounds1})}
If $\th_0\in\Theta^\b(M)$, then $h_n(\a;\th_0)$ is bounded above by
%
%
\begin{equation}
\label{EQUBbyINT} \qquad\frac{M}{N_\a\log N_\a} \sum_{i=1}^\infty
\frac{n^2i^{2\a-2\b
}(\log i)}{(i^{1+2\a+2p}+n)^2} \le\frac{M9e^{2\Aupper+3}}{N_\a\log N_\a
} \int_1^{N_\a}
x^{2\a
-2\b}\log x \,dx,
\end{equation}
by Lemma~\ref{LemmaTechLem10new} with $l=2$, $m=1$,
$s=2\a-2\b$ and, hence, $c=lr-s-1=1+2\a+2\b+4p\ge1$, for $n\ge
e^{4+8\a+8p}$.
Because the integrand is increasing in~$\a$, we find that
\begin{eqnarray*}
\sup_{\a\le\b-K/\log n}h_n(\a;\th_0) &\le&
M9e^{2\Aupper+3} \sup_{\a\le\b-K/\log n} N_\a^{-1}
\int_1^{N_\a} x^{-2K/\log n} \,dx
\\
&\le& M9e^{2\Aupper+3}\sup_{\a\le\b-K/\log n} \frac{N_\a
^{-2K/\log n}}{1-2K/\log n}
\\
&\le& M18e^{2\Aupper+3} e^{-2K/(1+2\b-2K/\log n+2p)},
\end{eqnarray*}
for $2K/\log n\le1/2$.
By\vspace*{1pt} its definition, $\underline\a_n\ge\b-K/\log n$
if the left-hand side of the preceding display is bounded above by $1/(16C^8)$.
This is true for $K\ge K_1$ as given in (\ref{EqKs}).
\end{pf*}

\begin{pf*}{Proof of (\ref{eqBoundsForBounds2})}
By Lemma~\ref{LemLowerBoundhn} (below), for $\th_0\in\Theta_{ss}^\b
(M,\eps)$
and $n\ge N_0^{1+2\Aupper+2p}$,
\begin{eqnarray*}
&&\inf_{\b+K/\log n\le\a\le\Aupper}h_n(\a;\th_0)
\\
&&\qquad \ge\inf_{\b+K/\log n\le\a\le\Aupper}\frac{\eps M}{(\rho
^{1+2\a+2p}+1)^2} n^{(2\a-2\b)/(1+2\a+2p)}
\\
&&\qquad \ge\frac{\eps M}{(\rho^{1+2\Aupper+2p}+1)^2}e^{2K/(1+2\Aupper+2p)}.
\end{eqnarray*}
By its definition $\overline\a_n\le\b+K/\log n$ if the right-hand side is
greater than $8C^8$. This happens for large enough $K\ge K_{2}$ as indicated.
\end{pf*}

%
\begin{lemma}
\label{LemLowerBoundhn}
For $\th_0\in\Theta_{ss}^\b(M,\eps)$ and $n\ge N_0^{1+2\a+2p}\vee e^4$,
\begin{eqnarray*}
h_n(\a;\th_0)&\ge& n^{(2\a-2\b)/(1+2\a+2p)}
\frac{\eps M}{(\rho
^{1+2\a+2p}+1)^2}.
\end{eqnarray*}
\end{lemma}

\begin{pf}
The function $h_n(\a;\th_0)$ is always bounded below by
\[
\frac{1}{N_\a\log N_\a}\sum_{N_\a\le i\le\rho N_\a}\frac{n^2
i^{1+2\a}(\log i)\th_{0,i}^2}{
(i^{1+2\a+2p}+n)^2} \ge
\frac{N_\a^{2\a}}{(\rho^{1+2\a+2p}+1)^2}\sum_{N_\a\le i\le
\rho N_\a}\th_{0,i}^2.
\]
For $\th_0\in\Theta_{ss}^\b(M,\eps)$
we can apply the definition of self-similarity to bound the sum on the
far right
below by $\eps M N_\a^{-2\b}$, for $N_\a\ge N_0$.
\end{pf}

\subsection{Proof of Proposition~\texorpdfstring{\protect\ref{lemmaPropositionRectangleMinimxRateEB}}{3.8}}
The proposition is an immediate consequence of (\ref{eqUBr}), (\ref
{eqratesrewrite}) and
(\ref{eqBoundsForBounds1}).

\section{Proof of Theorem \texorpdfstring{\protect\ref{TheoremSobolevRadiusPolishedTail}}{3.13}}\label{SectionProofSobolev}
It follows from the proof of Lemma 2.1 of \citet{KSzVZ} that, for
$\theta_0\in S^{\beta}(M)$,
\[
h_n(\a;\th_0)\leq Mn^{-\vafrac{1\wedge2(\beta-\alpha)}{1+2\a
+2p}}.\label{eqUBhn}
\]
The right-hand side is strictly smaller than $1/(16C^8)$
for $\a\leq\beta-2K/\log n$ with $K$ satisfying
\[
e^{2K/(1+2\beta-2K/\log n+2p)}=16C^8M.
\]
By the definition of $\underline\a_n$, we conclude that $\underline
\a_n\geq\beta-K/\log n$.
The theorem follows by combining this with (\ref{eqUBr}).

\section{Proof of Theorem \texorpdfstring{\protect\ref{TheoremSupersmoothRadiusPolishedTail}}{3.14}}
\mbox{}

\begin{pf*}{Proof of (\ref{rateC00})}
For $\a\leq\sqrt{\log n}/(3\sqrt{\log N_0})$ and
$\sqrt{\log n}\geq(3\sqrt{\log N_0})(1+2p)$, we have that
\[
N_\a\ge e^{\afrac{\log n}{1+2\sqrt{\log n}/(3\sqrt{\log
N_0})+2p}}\geq e^{\sqrt{\log N_0}\sqrt{\log n}}.
\]
Since also $N_0^{2\a}\leq e^{(2/3)\sqrt{\log N_0}\sqrt{\log n}}$,
the function $h_n(\a;\th_0)$ is bounded above by, for $\th_0\in
C^{00}(N_0,M)$,
\[
(\log N_0)N_0^{2+2\a}M/(N_\a
\log N_\a)\leq(\log N_0) N_0^2
M e^{-(1/3)\sqrt{\log N_0}\sqrt{\log n}}.
\]
Since this tends to zero, it will be smaller than $1/(16C^{8})$,
for large enough $n$, whence
$\underline\a_n\geq\sqrt{\log n}/(3\sqrt{\log N_0})$, by its definition.
Assertion (\ref{rateC00}) follows by substituting this
into (\ref{eqUBr}).
\end{pf*}

\begin{pf*}{Proof of (\ref{rateSinfty})}
As before, we give an upper bound for $h_n(\a;\th_0)$
by splitting the sum in its definition into two parts. The sum over the
indices $i > N_\a$ is bounded by
\[
\frac{1}{N_{\a}\log N_{\a}} n^2\sum_{i> N_\a}
i^{-1-2\a-4p}(\log i) \th_{0,i}^2.
\]
Since the function $f(x)=x^{-1-2\a-4p}\log x$ is monotonely decreasing
for $x\geq e^{1/(1+2\a+4p)}$, we have $N_\a^{-1}i^{-1-2\a-4p}(\log i)
\le N_\a^{-2-2\a-4p}(\log N_\a)= n^2N_\a^{2\a}\times\break \log N_\a$, for
$i>N_\a$.
Hence, the right-hand side is bounded by a multiple of
\[
N_\a^{2\a}\sum_{i>N_\a}^{\infty}
\th_{0,i}^2\le ne^{-cN_\a^d} \sum
_{i>N_\a}^{\infty}e^{cN_\a^d}\th_{0,i}^2
\le ne^{-cN_\a^d} M.
\]
The sum over the indices $i \leq n^{1/(1+2\a+ 2p)}$ is bounded by
\[
N_{\a}^{-1} \sum_{i \le N_\a}i^{1+2\a}
\th_{0,i}^2 \le N_{\a}^{-1}
e^{\vfrac{1+2\a}{d}\log\vafrac{1+2\a}{cd}}M,
\]
since the maximum on $(0,\infty)$ of the function $x \mapsto x^{1+2\a
}\exp(-c x^d)$ equals
$((1+2\a)/(cd))^{(1+2\a)/d}$.
Combining the two bounds, we find that for $\a\leq\sqrt{\log n}/\log
\log n$ and
sufficiently large $n$ the function $h_n(\a;\th_0)$ is bounded from
above by a multiple of
\[
M n\exp\bigl(-ce^{(d/3)\sqrt{\log n}\log\log n} \bigr)+ M \exp\biggl
(\frac{3\sqrt{\log n}}{2d}-
\frac{\sqrt{\log n}\log\log n}{3} \biggr).
\]
Since this tends to zero, the inequality $h_n(\a;\th_0) < 1/(16C^8)$
holds for large enough $n$ (depending on $p,C,c,d$ and $M$),
whence $\underline\a_n\ge\sqrt{\log n}/\log\log n$. Combining with
(\ref{eqUBr}),
this proves (\ref{rateSinfty}).
\end{pf*}

\section{Proof of Theorem~\texorpdfstring{\protect\ref{LemmaCounterExample}}{3.1}}
\label{secCounterExample}

From the proof of Theorem~\ref{thmmle} it can be seen that
the lower bound in (\ref{eqm1}) 
is valid also for non self-similar $\th_0$:
\begin{eqnarray*}
\inf_{\th_0\in\ell_2}\Pr\bigl(\underline\a_n(
\th_0)&\le&\hat\a_n \bigr)\rightarrow1. 
\end{eqnarray*}
In terms of the notation (\ref{EqDefWB}), we have that
$\th_0\in\hat C_n(L)$ if and only if $\llVert \hat\th_{\hat\a}-\th
_0\rrVert
\le L r(\hat\a)$,
which implies that $\llVert B(\hat\a;\th_0)\rrVert \le Lr(\hat\a
)+\llVert W(\hat\a)\rrVert $.
Combined with the preceding display, it follows that $P_{\th_0}
(\th_0\in\hat C_n(L) )$
is bounded above by
%
%
\begin{equation}
\qquad P_{\th_0} \Bigl(\inf_{\a\ge\underline\a_n(\th_0)}\bigl\llVert B(\a,
\th_0)\bigr\rrVert\leq L \sup_{\a\ge\underline\a_n(\th_0)}r(\a)+\sup
_{\a\ge\underline
\a_n(\th_0)}\bigl\llVert W(\a)\bigr\rrVert\Bigr) + o(1).
\label{eqCounterUBcoverage}
\end{equation}
The proofs of (\ref{w}) and (\ref{eqUBr}) show also that
%
%
\begin{eqnarray}
\label{eqrr} \sup_{\th_0\in\Theta^\b(M)}\sup_{\a\ge\underline\a_n(\th
_0)}r(\a) &\le& C_\b n^{-\afrac{\underline\a_n(\th_0)}{1+2\underline\a_n(\th_0)+2p}},
\\
\label{eqvv} \inf_{\th_0 \in\Theta^\b(M)}P \Bigl(\sup_{\a\ge\underline\a
_n(\th_0)}
\bigl\llVert W(\a)\bigr\rrVert&\le& C_3n^{-\afrac{\underline\a_n(\th
_0)}{1+2\underline\a_n(\th
_0) + 2p}} \Bigr)
\rightarrow1.
\end{eqnarray}
We show below that $\underline\a_{n_j}(\th_0)\ge\a_j$, for the
parameter $\th_0$ given in the theorem and
$\a_j$ the solution of
\[
\rho_j^{(1+2\beta+2p)/4}=n_j^{\a_j-\beta}.
\]
Since $\rho_j\rightarrow\infty$, it is clear that $\a_j>\beta$.
Furthermore, from
the assumption that $n_j\ge(2\rho_j^2)^{1+2\b+2p}n_{j-1}$, where
$n_{j-1}\ge1$,
it can be seen that $\a_j<\beta+1/2$.
By combining this with (\ref{eqrr}) and (\ref{eqvv}),
we see that the latter implies that the
expression to the right-hand side of the inequality in (\ref{eqCounterUBcoverage})
at $n=n_j$ is bounded above by a constant times
\begin{eqnarray*}
L n_j^{-\afrac{\a_j}{1+2\a_j + 2p}} &=& Ln_j^{-\afrac{(1+2p)(\a_j-\beta
)}{(1+2\beta+2p)(1+2\a_j+2p)}}n_j^{-\afrac{\beta}{1+2\beta+2p}}
\\
&\ll& L \rho_j^{-\afrac{(1+2p)}{4(2+2\beta+2p)}} n_j^{-\afrac{\beta
}{1+2\beta+2p}}.
\end{eqnarray*}
Since $\a\mapsto\llVert B(\a;\th_0)\rrVert $ is increasing, the
infimum on the
left-hand side of the
inequality is bounded from below by $\llVert B(\a_j;\th_0)\rrVert $.
Now, in view of (\ref{EqBiasSquare}), with $N_{j}=n_j^{1/(1+2\b+2p)}$
(and an abuse of notation,
as $N_\a$ was used with a different meaning before),
\begin{eqnarray*}
\bigl\llVert B(\a;\th_0)\bigr\rrVert^2&\ge&
\frac{1}{C^8} \sum_{N_j \le i <2N_j}\frac{i^{2+4\a+4p}\th
_{0,i}^2}{(i^{1+2\a+ 2p}+n_j)^2} \ge
\frac{1}{4C^8}\sum_{N_j\le i <2N_j}\th_{0,i}^2,
\end{eqnarray*}
for $\a\ge\b$, since $n_j \le i^{1+2\a+2p}$ for $i \ge N_j$. Using the
definition of $\th_0$, we see that this is lower bounded by a multiple
of
\[
M N_jn_j^{-(1+2\b)/(1+2\b+2p)} =M n_j^{-2\b/(1+2\b+2p)}.
\]
Thus, we deduce that the expression to the left of the
inequality sign in (\ref{eqCounterUBcoverage})
is of larger order than the expression to the right, whence the probability
tends to zero along the subsequence $n_j$.

Finally, we prove the claim that $\underline\a_{n_j}(\th_0)\ge\a
_j$, by showing that
$h_{n_j}(\a;\th_0)<1/(16C^8)$ for all $\a< \a_j$. We consider the cases
$\a\le\b$ and $\a\in(\b,\a_j)$ separately.
For $\a\leq\beta$ we have, by Lemma~\ref{LEMTAIL},
\begin{eqnarray*}
h_{n}(\a,\theta_0)&\leq&\frac{1}{N_\a\log N_\a} \Biggl(\sum
_{i=1}^{N_\a}M\log i +n^2\sum
_{i\geq N_\a}Mi^{-2-2\a-4p-2\beta}\log i \Biggr)\lesssim M.
\end{eqnarray*}
Thus, $h_n(\a,\th_0)$ is smaller than $1/(16C^8)$ for sufficiently
small $M$.
For $\beta<\a<\a_j$ we have that
$h_{n_j}(\a;\th_0) \le A_1+A_2+A_3$ for
\begin{eqnarray*}
A_1&=&\frac{1+2\a+2p}{n_j^{1/(1+2\a+2p)}\log n_{j}}\sum_{i \le\rho
_j^{-1}N_j}Mi^{2\a-2\b}
\log i,
\\
A_2&=& \frac{1+2\a+2p}{n_j^{1/(1+2\a+2p)}\log n_{j}} \sum_{N_j\le i < 2N_{j}}
\frac{n_j^2i^{1+2\a}(\log i)N_j^{-1-2\b
}M}{(i^{1+2\a+2p}+n_j)^2},
\\
A_3&=&\frac{1+2\a+2p}{n_j^{1/(1+2\a+2p)}\log n_{j}}\sum_{i \ge\rho
_jN_j}Mn_{j}^2
i^{-2-2\a-4p-2\b}(\log i).
\end{eqnarray*}
The first term satisfies
\begin{eqnarray*}
A_1&\lesssim& M\rho_j^{-(1+2\a-2\beta)}n_j^{\vafrac{1+2\a-2\beta
}{1+2\beta+2p}-\afrac{1}{1+2\a+2p}}
\\
&\le& M\rho_j^{-1}n_j^{\afrac{2(\a-\beta)(2+2\a+2p)}{(1+2\beta
+2p)(1+2\a+2p)}}.
\end{eqnarray*}
The right-hand side is increasing in $\a$, and hence is maximal over $(\b
,\a_j]$ at $\a_j$.
At this value it tends to zero in view of the definition of $\a_j$.
By Lemma~\ref{LEMTAIL},
\begin{eqnarray*}
A_3&\lesssim& M\rho_j^{-(1+2\a+4p+2\beta)}n_j^{2-\afrac{1}{1+2\a
+2p}}N_j^{-1-2\a-4p-2\b}
\\
&\le& M\rho_j^{-1}n_j^{\afrac{2(\b-\a)(2\a+2p)}{(1+2\beta+2p)(1+2\a+2p)}}.
\end{eqnarray*}
This tends (easily) to zero for $\a>\b$.\vspace*{1pt}

The term $i^{1+2\a+2p}+n_j$ in the denominator of the sum in $A_2$ can
be bounded
below both by $i^{1+2\a+2p}$ and by $n_j$, and there are at most $N_j$
terms in the sum.
This shows that
\begin{eqnarray*}
A_2 &\lesssim&\frac{n_j^{-1/(1+2\a+2p)}}{\log n_{j}} N_j \biggl(
\frac{n_j^2}{N_j^{1+2\a+4p}}\wedge(2N_j)^{1+2\a} \biggr)
\log(2N_j)N_j^{-1-2\b}M
\\
& \lesssim&  M \bigl(n_j^{\vafrac{1+2\b- 2\a}{1+2\b+ 2p}- \afrac{1}{1+2\a+
2p}}
\\
&&\hspace*{13pt}{} \wedge n_j^{\vafrac{1+2\a- 2\b}{1+2\b+ 2p}- \afrac{1}{1+2\a+
2p}}
\bigr).
\end{eqnarray*}
The exponents of $n_j$ in both terms in the minimum are equal to 0 at
$\a=\b$.
For $\a\ge\b$ the first exponent is negative, whereas the second
exponent is increasing
in $\a$ and hence negative for $\a<\b$. It follows that $A_2
\lesssim M$.

Putting things together, we see that
$\limsup_{j\rightarrow\infty} \sup_{\a\le\a_j}h_{n_j}(\a; \th_0)$
can be made arbitrarily small by choosing $M$ sufficiently small.

\section{Technical lemmas}
\label{SectionTechnicalLemmas}

The following two lemmas are Lemma 8.3 and an extension of Lemma 8.2 of
\citet{KSzVZ}.

%
\begin{lemma}\label{LemmaTechLem6}
For any $p \geq0$, $r \in(1, (\log n)/(2\log(3e/2))]$, and $g > 0$,
\[
\sum_{i=1}^{\infty}\frac{n^g\log i}{(i^{r}+n)^g} \geq
\frac
{1}{3\cdot2^{g}}n^{1/r}(\log n/ r).
\]
\end{lemma}

%
\begin{lemma}\label{LemmaTechLem10}
For any $l, m, r, s\ge0$ with $c:=lr-s-1>0$ and $n \ge e^{(2mr/c)\vee r}$,
\[
\bigl(3^r+1\bigr)^{-l} ({\log n}/r )^{m}n^{-c/r}
\le\sum_{i=1}^\infty\frac{i^{s}(\log i)^m}{(i^{r}+n)^l} \le
\bigl(3+2c^{-1}\bigr) ({\log n}/r )^{m}n^{-c/r}.
\]
\end{lemma}

The series in the preceding lemma changes in character as
$s$ decreases to $-1$, with the transition starting at $-1+O(1/\log n)$.
This situation plays a role in the proof of the key result
(\ref{eqDiffBounds}) and is handled by the following lemma.

%
\begin{lemma}\label{LemmaTechLem10new}
For any $l, m, r\ge0$ and $s\in\RR$ with $c:=lr-s-1>0$ and $n \ge
e^{(2mr/c)\vee(4r)}$,
\[
\sum_{i=2}^\infty\frac{i^{s}(\log i)^m}{(i^{r}+n)^l} \le
e^{\llvert s\rrvert +m+2} 3\bigl(1+2c^{-1}\bigr) n^{-l} \int
_1^{n^{1/r}} x^s(\log x)^m
\,dx.
\]
\end{lemma}

\begin{pf}
For $N=n^{1/r}$ the series is bounded above by $I+\mathit{II}$, for
\[
I=n^{-l}\sum_{i\le N} i^s(\log
i)^m,\qquad \mathit{II}= \sum_{i> N}
i^{s-rl} (\log i)^m.
\]
We treat the two terms separately.

Because the function $f\dvtx  x\mapsto x^s (\log x)^m$ is monotone or
unimodal on $[2,N]$,
the sum $I$ is bounded by $2f(\mu)+\int_2^N f(x) \,dx$, for $\mu$ the
point of maximum. As
the derivative of $\log f$ is bounded above by $\llvert s\rrvert +m$,
it follows that
$f(x)\ge e^{-\llvert s\rrvert -m} f(\mu)$
in an interval of length 1 around the point of maximum, and hence
$f(\mu)$ is bounded by
$e^{\llvert s\rrvert +m}$ times the integral.

By Lemma~\ref{LEMTAIL}, with $k=rl-s-1=c$, the term $\mathit{II}$ is bounded above
by $(1+2c^{-1})(\log N)^mN^{-c}$, for $N \ge e^{2m/c}$. This is bounded
by the right-hand side of the lemma if
\[
(\log N)^m N^{s+1}\le\bigl(1+\llvert s\rrvert\bigr)
e^{m+2}\int_1^N x^s(
\log x)^m \,dx.
\]
For $\llvert s+1\rrvert \le1/\log N$, we have $N^{s+1}\le e$ and
$x^s/x^{-1}\ge
e^{-1}$ for $x\in[1,N]$.
The former bounds the left-hand side by $(\log N)^m e$, while the latter gives
that the integral on the right is bounded below by $e^{-1}\int_1^N
x^{-1}(\log x)^m \,dx=
e^{-1}(m+1)^{-1}(\log N)^{m+1}$, whence the display is valid.
For $s+1\ge1/\log N$, we bound the integral below by the
integral of the same function over $[\sqrt N,N]$, which is bounded
below by $(\log\sqrt N)^m(N^{s+1}-N^{(s+1)/2})/(s+1)\ge(\log\sqrt
N)^m N^{s+1}/ (4(s+1))$,
as $(3/4) N^{s+1}\ge N^{(s+1)/2}$ if $s+1\ge1/\log N$. As also
$(1+\llvert s\rrvert )/(1+s)\ge1$, this
proves the display. For $s+1\le-1/\log N$, we similarly bound the
integral below by
$(\log\sqrt N)^m(N^{(s+1)/2}-N^{s+1})/\llvert s+1\rrvert \ge(\log
\sqrt N)^m
N^{s+1}/(4\llvert s+1\rrvert )$.
\end{pf}

%
\begin{lemma}
\label{LEMTAIL}
For $k>0$, $m\ge0$ and $N \ge e^{2m/k}$,
\[
\sum_{i>N}i^{-1-k}(\log i)^m
\le(1/N+2/k) (\log N)^m N^{-k}.
\]
\end{lemma}

\begin{pf}
Because the function $x\mapsto x^{-1}\log x$ is decreasing for $x\ge
e$, we have
$i^{-k/2}(\log i)^m\le N^{-k/2}(\log N)^m$ for
$i^{k/(2m)}>N^{k/(2m)}\ge e$ if $m>0$. If $m=0$, this
inequality is true for every $i> N\ge1$. Consequently,\vspace*{1pt} the sum in the
lemma is bounded above by
$N^{-k/2}(\log N)^m\sum_{i>N}i^{-1-k/2}$. The last\vspace*{1pt} sum is bounded
above by $N^{-1-k/2}
+\int_N^\infty x^{-1-k/2} \,dx=(1/N+2/k)N^{-k/2}$.
\end{pf}

\section{Concluding remarks}
\label{seccon}
A full Bayesian approach, with a hyperprior on $\a$, is an alternative
to the empirical Bayesian approach employed here. As the full posterior
distribution
is a mixture of Gaussians with different means, there are multiple
reasonable definitions
for a credible set. Based on our earlier work on rates of contraction
[\citet{KSzVZ}], we believe that their coverage
will be similar to the empirical Bayes sets considered in the present paper.

Rather than balls, one may, in both approaches, consider sets of
different shapes,
for instance, bands if the parameter can be identified with a function.
It has already been
noted in the literature that rates of contraction of functionals, such
as a function
at a point, are suboptimal unless the prior is made dependent on the functional.
Preliminary work suggests that adaptation complicates this situation further,
except perhaps when the parameters are self-similar.

The question of whether a restriction to polished tail or self-similar
sequences is reasonable from
a practical point of view is open to discussion. From the theory and
the examples
in this paper it is clear that a naive
or automated (e.g., adaptive) approach will go wrong in certain situations.
This appears to be a fundamental weakness of statistical uncertainty
quantification:
honest uncertainty quantification is always conditional on a set of assumptions.
To assume that the true sequence is of a polished tail type is reasonable,
but it is not obvious how one would communicate this assumption to a
data analyst.

\begin{appendix}
\section*{Appendix: Proof of (\texorpdfstring{\protect{\lowercase{\ref{eqm1}}}}{5.4}) in Theorem~\texorpdfstring{\protect{\lowercase{\ref{thmmle}}}}{5.1}}
\label{secproofMLE2}
With the help of the dominated convergence theorem
one can see that the random function $\ell_n$ is differentiable and
the derivative is given by
%
%
\begin{equation}
\label{eqm} \MM_n(\a) = \sum_{i=1}^{\infty}
\frac{n\log i}{i^{1+2\a}\k_i^{-2}+n} -\sum_{i=1}^{\infty}
\frac{n^2i^{1+2\a}\k_i^{-2}\log i}{(i^{1+2\a
}\k_i^{-2}+n)^2}X_i^2.
\end{equation}
The proof of (\ref{eqm1}) consists of the following steps:
\begin{longlist}[(ii)]
\item[(i)] In Section~\ref{secLB} we show that with probability
tending to one, uniformly over
$\th_0\in\ell_2$,
the function $\MM_n$ is strictly positive on the interval
$(0,\underline\a_n)$.
\item[(ii)] In Section~\ref{secUB} we show that on the interval
$[\overline\a_n,\Aupper)$ the process $\MM_n$ is strictly negative
with probability
tending to one, uniformly over $\th_0\in\Theta_{pt}(L_0)$.
\end{longlist}
Steps (i) and (ii) show that $\MM_n$ has no local maximum
on the respective interval.

\subsection{The process $\mathbb{M}_n$ on \texorpdfstring{$(0,\underline\alpha_n]$}{(0,underlinealphan]}}
\label{secLB}
We can assume $\underline\a_n>0$, which leads to the
inequality $h_n(\a;\th_0)\leq(16C^8)^{-1}$ for every $\a\in
(0,\underline\a_n]$.
For the proof of (i) above it is sufficient to show that the following hold:
%
%
\begin{eqnarray}
\liminf_{n\rightarrow\infty}\inf_{\th_0\in\ell_2}\inf
_{\a\in
(0,\underline\a_n]} E_{\th_0}\frac{\MM_n(\a)(1+2\a+2p)}{n^{\afrac{1}{1+2\a+2p}}\log
n}&>&
\frac{1}{48 C^4},\label{eqMain3}
\\[-2pt]
\sup_{\th_0\in\ell_2}E_{\th_0}\sup_{\a\in(0,\underline\a
_n]}
\frac{\llvert \MM_n(\a)-\EE_0
\MM_n(\a)\rrvert (1+2\a+2p)}{n^{\afrac{1}{1+2\a+2p}}\log n}&\rightarrow&
0.\label{eqMain4}
\end{eqnarray}
The expectation in (\ref{eqMain3}) is equal to
\begin{eqnarray*}
&& \frac{1+2\a+2p}{n^{\afrac{1}{1+2\a+2p}}\log n} \Biggl(\sum_{i=1}^{\infty}
\frac{n^2\log i}{(i^{1+2\a}\k_i^{-2}+n)^2} -\sum_{i=1}^{\infty}
\frac{n^2i^{1+2\a}(\log i)\th
_{0,i}^2}{(i^{1+2\a}\k_i^{-2}+n)^2} \Biggr)
\\[-2pt]
&&\qquad \geq\frac{1+2\a+2p}{C^4n^{\afrac{1}{1+2\a+2p}}\log n}\sum_{i=1}^{\infty}
\frac{n^2\log
i}{(i^{1+2\a+2p}+n)^2} -C^4h_n(\a;\th_0).
\end{eqnarray*}
By Lemma~\ref{LemmaTechLem6} [with $g=2$, $r=1+2\a+2p$ and $\log
n\geq(8\log(3e/2))^2\vee4(1+2p)\log(3e/2)$],
the first term of the preceding display is bounded from below
by $1/(12C^4)$ for all $\a\in(0,\underline\a_n)\subset(0,\sqrt
{\log n}/4)$.
Inequality (\ref{eqMain3}) follows, as
the second term is bounded above by $C^4/(16C^8)$, by the definition of
$\underline\a_n$.

To verify (\ref{eqMain4}), it (certainly) suffices
by Corollary~2.2.5 in \citet{vdVW} applied with $\psi(x)=x^2$
to show that for any positive $\a\leq\underline\a_n\leq\Aupper$
%
%
\begin{eqnarray}
\Var_{\th_0} \frac{\MM_n(\a)(1+2\a+2p)}{n^{1/(1+2\a+2p)}\log n} &\leq&
K_1e^{-(3/2)\sqrt{\log n}}\label{eqVarMn3},
\\
\int_0^{\diam} \sqrt{N\bigl(\eps, (0,\underline
\a_n],d_n\bigr)} \,d\eps&\leq& K_2
e^{-(9/8)\sqrt{\log n}}(\log n),\label{eqIntCovNumb}
\end{eqnarray}
where $d_n$ is the semimetric defined by
\[
d_n^2(\a_1, \a_2) =
\Var_{\th_0} \biggl(\frac{\MM_n(\a_1)(1+2\a
_1+2p)}{n^{1/(1+2\a_1+2p)}\log n} -\frac{\MM_n(\a_2)(1+2\a
_2+2p)}{n^{1/(1+2\a_2+2p)}\log n} \biggr),
\]
$\diam$ is the diameter of $(0,\underline\a_n]$ relative to $d_n$,
$N(\eps,B,d_n)$ is the minimal number of
$d_n$-balls of radius $\eps$ needed to cover the set $B$, and the
constants $K_1$ and $K_2$ do not
depend on the choice of $\th_0\in\ell_2$.

By Lemma 5.2 of \citet{KSzVZ}, the variance in (\ref{eqVarMn3}) is bounded
above by a multiple of, for any $\a\in(0,\underline\a_n)\subset
(0,\sqrt{\log n}/4)$,
\begin{eqnarray*}
n^{-1/(1+2\a+2p)} \bigl(1+h_n(\a;\th_0) \bigr) &\le&
n^{-1/(1+2\sqrt{\log n}/4+2p)} \bigl(1+\bigl(16C^8\bigr)^{-1} \bigr)
\\
&\le&\bigl(1+\bigl(16C^8\bigr)^{-1} \bigr)
e^{-(3/2)\sqrt{\log n}},
\end{eqnarray*}
%
for $\log n\geq(6(1+2p))^2$. Combination with the triangle inequality shows
that the $d_n$-diameter of the set
$(0,\underline\a_n)$ is bounded by a constant times $e^{-(3/4)\sqrt
{\log n}}$.
To verify (\ref{eqIntCovNumb}), we apply Lemma 5.3 of \citet{KSzVZ}
according to which,
for any $0<\a_1<\a_2<\infty$,
\begin{eqnarray*}
&& \Var_{\th_0} \biggl( \frac{(1+2\a_1+2p)\MM_n(\a_1)}{n^{1/(1+2\a
_1+2p)}\log n} -\frac{(1+2\a_2+2p)\MM_n(\a_2)}{n^{1/(1+2\a_2+2p)}\log
n} \biggr)
\\
&&\qquad \lesssim(\a_1-\a_2)^2(\log
n)^2\sup_{\a\in[\a
_1,\a_2]}n^{-1/(1+2\a+2p)}
\bigl(1+h_n(\a;\th_0) \bigr).
\end{eqnarray*}
%
We see that for $\a_1,\a_2\in(0,\underline\a_n]$ the metric
$d_n(\a_1,\a_2)$ is bounded above by a constant times
$(\log n) e^{-(3/4)\sqrt{\log n}}\llvert \a_1-\a_2\rrvert $.
Therefore, the coverage\break  number of the interval $(0,\underline\a_n)$
is bounded above by a constant times\break
$(e^{-(3/4)\sqrt{\log n}}(\log n)^{3/2})/\eps$, which leads to the inequality
\[
\int_0^\diam\sqrt{N\bigl(\eps, (0,\underline
\a_n], d_n\bigr)} \,d\eps\lesssim e^{-(9/8)\sqrt{\log n}}(\log
n)^{3/4},
\]
where the multiplicative constant does not depend on the choice of $\th_0$.

\subsection{The process $\MM_n$ on \texorpdfstring{$(\overline\a_n,\Aupper]$}
{(overlinealphan,A]}}\label{secUB}

To prove that $\ell_n$ is strictly decreasing on $(\overline\a
_n,\Aupper]$,
it is sufficient to verify the following:
%
%
\begin{eqnarray}
\qquad \limsup_{n\rightarrow\infty}\sup_{\th_0\in\Theta_{pt}(L_0)}\sup
_{\a\in(\overline\a_n,\Aupper]} E_{\th_0}\frac{\MM_n(\a)(1+2\a
+2p)}{n^{\afrac{1}{1+2\a+2p}}h_n(\a
;\th_0)\log n }&<&-
\frac{3}{8C^4}, \label{eqMain1}
\\
\sup_{\th_0\in\Theta_{pt}(L_0)}E_{\th_0}\sup_{\a\in(\overline
\a_n,\Aupper]}
\frac{\llvert \MM_n(\a)-
E_{\th_0} \MM_n(\a)\rrvert (1+2\a+2p)}{n^{\afrac{1}{1+2\a+2p}}h_n(\a;\th
_0)\log n }&\rightarrow&0. 
\label{eqMain2}
\end{eqnarray}
We shall verify this under the assumption that $\overline\a_n<
\Aupper$,
using that in this case $h_n(\a;\th_0)\geq8C^8$, for all $\a\in
[\overline\a_n,\Aupper]$, by the definition of $\overline\a_n$.

In view of (\ref{assumpinverse}), the expectation in (\ref{eqMain1})
is bounded above by
\[
\frac{(1+2\a+2p) C^4}{n^{\afrac{1}{1+2\a+2p}}8C^8\log n} \sum
_{i=1}^{\infty}
\frac{n^2\log i}{(i^{1+2\a+2p}+n)^2}-C^{-4}. 
\]
Inequality (\ref{eqMain1}) follows by an application of
Lemma~\ref{LemmaTechLem10} (with $s=0$, $r=1+2\a+2p$, $l=2$, $m=1$, and
hence $c=1+4\a+4p\ge1$).

To verify (\ref{eqMain2}), it suffices,
by Corollary~2.2.5 in \citet{vdVW} applied with $\psi(x)=x^2$,
to show that
%
%
\begin{eqnarray}
\Var_{\th_0} \frac{(1+2\a+2p)\MM_n(\a)}{n^{\afrac{1}{1+2\a
+2p}}h_n(\a;\th_0)\log n} &\leq& K_1
e^{-(3/2)\sqrt{\log n}},\label{eqVarMn1}
\nonumber\\[-8pt]\\[-8pt]\nonumber
\int_0^{\diam} \sqrt{N\bigl(\eps, [\overline
\a_n,\Aupper],d_n\bigr)} \,d\eps&\leq& K_2(
\log n)^{5/4}L_0^{1/2} e^{-(7/8)\sqrt{\log n}},
\nonumber
\end{eqnarray}
where this time $d_n$ is the semimetric defined by
\begin{eqnarray*}
&& d_n^2(\a_1, \a_2)
\\
&&\qquad =
\Var_{\th_0} \biggl(\frac{\MM_n(\a_1)(1+2\a
_1+2p)}{n^{\afrac{1}{1+2\a_1+2p}}h_n(\a_1;\th_0)\log n} -\frac{\MM
_n(\a_2)(1+2\a_2+2p)}{n^{\afrac{1}{1+2\a_2+2p}}h_n(\a_2;\th_0)\log
n} \biggr),
\end{eqnarray*}
and the constants $K_1$ and $K_2$ do not depend on $\th_0\in\Theta
_{pt}(L_0)$.

By Lemma 5.2 of \citet{KSzVZ}, the variance (\ref{eqVarMn1}) is
bounded above by a multiple of
\[
n^{-1/(1+2\a+2p)} \bigl(1+h_n(\a;\th_0)
\bigr)/h_n(\a;\th_0)^2 \lesssim
e^{-(3/2)\sqrt{\log n}},
\]
for $\a\in[\overline\a_n,\Aupper]$, since $h_n(\a;\th_0)\geq
8C^8$ and $\Aupper\leq\sqrt{\log n}/4$.
Combination with the triangle inequality shows that the $d_n$-diameter
of the set
$[\overline\a_n,\Aupper)$ is of the square root of this order.

By Lemma~\ref{Lemmametric1} below, the present metric $d_n$ is
bounded above
similarly to the metric $d_n$ in Section~\ref{secLB}.
The entropy number of the interval $(0,\overline\a_n)\subset(0,\sqrt
{\log n})$ is bounded above by $L_0(\log n)^{5/2}e^{-(1/4)\sqrt{\log
n}}$. Therefore, the corresponding entropy integral can also be bounded
in a similar way by a multiple of $L_0^{1/2}(\log n)^{5/4}
e^{-(7/8)\sqrt{\log n}}$.

%
\begin{lemma}\label{Lemmametric1}
For any $\th_0\in\Theta_{pt}(L_0)$, and any $0<\overline\a_n \le
\a_1<\a_2\leq\Aupper$,
\begin{eqnarray*}
&& \Var_{\th_0} \biggl[\frac{\MM_n(\a_1)(1+2\a_1+2p)}{n^{\afrac{1}{1+2\a_1+2p}}h_n(\a_1;\th_0)} -\frac{\MM_n(\a_2)(1+2\a
_2+2p)}{n^{\afrac{1}{1+2\a_2+2p}}h_n(\a
_2;\th_0)} \biggr]
\\
&&\qquad \leq
K L_0^2\frac{\llvert \a_1-\a_2\rrvert ^2(\log n)^4}{e^{(1/2)\sqrt
{\log n}}},
\end{eqnarray*}
where the constant $K$ does not depend on $\th_0\in\Theta_{pt}(L_0)$.
\end{lemma}

\begin{pf}
The left-hand side of the lemma can be written
$n^4\sum_{i=1}^{\infty} (f_i(\a_1)-f_i(\a_2) )^2 \Var
_{\th_0} X_i^2$, for
%
%
\begin{equation}
\label{eqmetric1} f_i(\a)=\frac{(1+2\a+2p)}{n^{1/(1+2\a+2p)}}\frac
{i^{1+2\a}\k
_i^{-2}\log i}{
(i^{1+2\a}\k_i^{-2}+n)^{2}h_n(\a;\th_0)}.
\end{equation}
The absolute value of the derivative $f_i'(\a)$ is equal to
\begin{eqnarray*}
&& f_i(\a)\biggl\llvert\frac{2}{1+2\a+2p}+\frac{2\log n}{(1+2\a
+2p)^2}
\\
&&\hspace*{27pt}{} +2\log
i-4\frac{(\log i)
i^{1+2\a}\k_i^{-2}}{
i^{1+2\a}\k_i^{-2}+n} -\frac{h'_n(\a;\th_0)}{h_n(\a;\th_0)}\biggr\rrvert
\\
&&\qquad \lesssim L_0\rho^{1+2\a}f_i(\a) (\log
i +\log n),
\end{eqnarray*}
by Lemma~\ref{LemmaDerivative}.
Writing the difference $f_i(\a_1)-f_i(\a_2)$ as the integral of
$f_i(\a)$, applying the
Cauchy--Schwarz inequality to its\vspace*{1pt} square, interchanging the sum and integral,
and substituting $\Var_{\th_0}X_i^2=2/n^2+4\k_i^2\th_{0,i}^2/n$,
we can bound the variance in the lemma by a multiple of
\[
(\a_1-\a_2)^2n^4L_0^2
\rho^{2+4\a_2} \sup_{\a\in[\a_1, \a_2]} \sum
_{i=1}^\infty f_i(\a)^2(\log
i+\log n)^2 \biggl(\frac{2}{n^2}+\frac
{4\k_i^2\th_{0,i}^2}{n} \biggr).
\]
The series splits in two terms by the last plus sign on the right.
Using Lemma~\ref{LemmaTechLem10}
with $s=2+4\a+4p$, $l=4$, $r=1+2\a+2p$ and $m=4$ or $m=2$
on\vspace*{1pt} the first part, and the inequality $ni^r(r\log i)^m/(i^r+n)^2\le
(\log n)^m$ for $n\ge e^4$,
with $r=1+2\a+2p$ and $m=3$ and $m=1$ on the second part, we
can bound the preceding display by a multiple of
\[
(\log n)^4 L_0^2\rho^{2+4\a_2}\sup
_{\a\in[\a_1,\a
_2]}n^{-1/(1+2\a+2p)} \bigl(1+h_n(\a;
\th_0) \bigr)/h_n(\a;\th_0)^2.
\]
We complete the proof by using that $h_n(\a;\th_0)\geq8C^8$ for
$\sqrt{\log n}/(4\sqrt{\log\rho})\geq\Aupper\geq\a\geq\a
_1\geq\overline\a_n$.
\end{pf}

%
\begin{lemma}
\label{LemmaDerivative}
For $\th_0\in\Theta_{pt}(L_0)$, $n\ge e^{(\log(N_0\rho)/3)^2}$ and
$\a\leq\sqrt{\log n}/4$,
\begin{eqnarray*}
h_n'(\a;\th_0)&\le& 48
\rho^{1+2\a} L_0(\log n) h_n(\a;
\th_0) .
\end{eqnarray*}
\end{lemma}

\begin{pf}
The derivative $h_n(\a;\th_0)$ can be computed to be
\[
\frac{2(1+\log N) h_n(\a;\th_0)}{1+2\a+2p} +\frac{1}{N\log N}\sum
_{i=1}^\infty
\frac{2n^2(\log i)^2(n-i^{1+2\a+2p})
i^{1+2\a}\th_{0,i}^2}{(i^{1+2\a+2p}+n)^3}.
\]
The series in the second term becomes bigger if we bound
$(n-i^{1+2\a+2p})/ \break (i^{1+2\a+2p}+n)$ by 1. Next, the series
can be split in the terms with $i\le N$ and $i>N$.
In the first one factor $\log i$ can be bounded by $\log N$, and hence
this part is bounded above by $(\log N) h_n(\a;\th_0)$. For $\th
_0\in\Theta_{pt}(L_0)$
the second part is bounded above by
\begin{eqnarray*}
\frac{2n^2}{N\log N}\sum_{i>N}(\log
i)^2i^{-1-2\a-4p} \theta_{0,i}^{2} &\le& 6
n^2N^{-2-2\a-4p}(\log N)\sum_{i>N}
\theta_{0,i}^{2}
\\
&\le& 6L_0(\log N)\rho^{1+2\a}\sum
_{N/\rho\leq i<N}i^{1+2\a}\theta_{0,i}^{2}
\\
&\leq& 48L_0 \rho^{1+2\a}(\log n)h_n(\a,
\theta_0),
\end{eqnarray*}
in view of Lemma~\ref{LEMTAIL} with $m=2$ and $k=1+2\a+4p\ge1$, for
$n\ge4$ and $N\geq e^{3\sqrt{\log n}}\geq\rho N_0$.
\end{pf}
\end{appendix}



%

\printaddresses
\end{document}